\documentclass[11pt,a4paper]{article}
\usepackage[utf8]{inputenc}
\usepackage{a4wide} 
\usepackage[utf8]{inputenc}
\usepackage{subfig}
\usepackage{bbm}
\usepackage{amsmath}
\usepackage{graphicx}
\usepackage{epstopdf}
\usepackage{algcompatible}
\usepackage{adjustbox}
\usepackage{tikz}
\usetikzlibrary{shapes,positioning,backgrounds}
\RequirePackage{ifthen}
\RequirePackage[T1]{fontenc}
\RequirePackage[ngerman,english]{babel}
\RequirePackage{amsmath}
\RequirePackage{amssymb}
\RequirePackage{amsthm}
\RequirePackage{amsfonts}
\RequirePackage{graphicx}
\RequirePackage{color}
\RequirePackage{listings} 
\usepackage{epstopdf}

\usepackage{algcompatible}
\RequirePackage[section,Algorithm]{algorithm}
\usepackage{a4wide} 
\usepackage[utf8]{inputenc}
\usepackage{subfig}
\usepackage{bbm}
\usepackage{amsmath}
\usepackage{graphicx}
\usepackage{epstopdf}
\usepackage{algcompatible}
\usepackage{adjustbox}
\usepackage{titlesec}
\usepackage{csquotes}
\RequirePackage{ifthen}
\RequirePackage[T1]{fontenc}
\RequirePackage[ngerman,english]{babel}
\RequirePackage{amsmath}
\RequirePackage{amssymb}
\RequirePackage{amsthm}
\RequirePackage{amsfonts}
\RequirePackage{graphicx}
\RequirePackage{color}
\RequirePackage{listings} 
\usepackage{epstopdf}

\usepackage{algcompatible}
\RequirePackage[section,Algorithm]{algorithm}
\RequirePackage{emptypage}

\usepackage{biblatex}
\renewbibmacro{in:}{}
\DeclareFieldFormat[article]{title}{\textit{#1}} 
\DeclareFieldFormat[article]{journaltitle}{#1}
\theoremstyle{plain}						
\newtheorem{theorem}{Theorem}[section]
\newtheorem{lemma}[theorem]{Lemma}

\newtheorem{definition}[theorem]{Definition}

\newtheorem{remark}[theorem]{Remark}
\numberwithin{equation}{section}

\titleformat{\section}[block]{\normalfont\bfseries}{\thesection.}{0.5em}{}
\titlespacing{\section}{0pc}{1pc}{1pc}

\titleformat{\subsection}[block]{\normalfont\bfseries}{\thesubsection.}{0.5em}{}
\titlespacing{\subsection}{0pc}{1pc}{1pc}

\begin{document}
	\title{Data-inspired modeling of accidents in traffic flow networks using the Hawkes process}
	\author{Simone G\"ottlich\footnotemark[1], \; Thomas Schillinger\footnotemark[1]}
	
	\footnotetext[1]{University of Mannheim, Department of Mathematics, 68131 Mannheim, Germany (goettlich@uni-mannheim.de, schillinger@uni-mannheim.de)}

	\date{ \today }
	
	\maketitle
	
	\begin{abstract}
		\noindent
		We consider hyperbolic partial differential equations (PDEs) for a dynamic description of the traffic behavior in road networks. These equations are coupled to a Hawkes process that models traffic accidents taking into account their self-excitation property which means that accidents are more likely in areas in which another accident just occurred. We discuss how both model components interact and influence each other. A data analysis reveals the self-excitation property of accidents and determines further parameters. Numerical simulations using risk measures underline and conclude the discussion of traffic accident effects in our model.
	\end{abstract}
	
	{\bf AMS Classification.} 35R60, 90B20, 65M06
	
	{\bf Keywords.} traffic flow network model, random accidents, Hawkes process, numerical simulations 
	
	
	
	\section{Introduction}
	Traffic accident modeling has been a frequently discussed topic in the last years because accidents are highly costly and lead to severe injuries or even death. There is a large variety of approaches that try to predict or understand accident occurrences in traffic flow, from statistical points of view \cite{Kwon2006,Le2020,Sipos2017} or queuing theory \cite{Jin2019}, where accidents are interpreted as a stochastic process. Another approach to describe accidents has been made by constructing Bayesian networks \cite{Mora2017,Zou2017}. Combining accident theory with traffic flow models, conservation laws have been exploited to describe collisions \cite{Moutari2014,Moutari2013}. A bi-directional relation between accident behaviour and traffic flow using conservation laws has been established first in \cite{Goettlich2020,Goettlich2021}.
	Unlike the already existing models, our approach is motivated by the self-excitation property of accidents, which means that further accidents are more likely to happen when there has just been an accident. A suitable tool to represent self-excitation is the Hawkes process which was first introduced by Alan G. Hawkes in 1971, see \cite{Hawkes1971}. It belongs to the class of counting processes and has the self-exciting property. This typically leads to a clustering in the jumps of the process that correspond to accident in our model. As a first application 20 years ago, the Hawkes process has been applied to model occurrences of earthquakes \cite{Ogata1999} where aftershocks are usual phenomena. There have been various further areas in which the Hawkes process has been applied, such as for the clustering of crimes \cite{Mohler2011}, failures in supply networks \cite{Verheugd2020} or several issues concerning financial markets, where for example credit defaults or optimal execution strategies have been modeled using the Hawkes process \cite{Bacry2015, Errais2010,Zagst2020} {or limit order book dynamics have been investigated \cite{Chen2017}}. Recently, the Hawkes process has also been used to describe the evolution of the Covid-19 pandemic \cite{Garetto2021}. 
	In the context of traffic flow the Hawkes process has been considered only in a few works so far. In \cite{Lim2016} the clustering effect of traffic flow and congestion was described using the Hawkes process. Further works proposed a traffic accident model accounting for multiple collisions \cite{Li2018}. Therein, also a brief parameter estimation on artificially generated accident data was performed. A real data study on a Hawkes process accident model in London and Rome including a parameter estimation is presented in \cite{diLoro.2024,Kalair2020}. Contrary, the work in \cite{Goettlich.2024} follows a probabilistic approach inspired by the Hawkes process including a self-excitation property but focuses on long-term accident trends.
	
	However, these works model traffic accidents independently from traffic flow. In this paper, we set up a rigorous mathematical model based on a hyperbolic conservation law describing the evolution of the traffic density. Accidents are introduced by a capacity reduction on roads leading to a space-dependent flux function in the conservation law. We consider two classes of accidents: On the one hand, accidents that are directly linked to the traffic flow, meaning that there is a larger accident risk when there is a high traffic density and high velocity. This type of accidents represent a background noise and have been originally introduced in \cite{Goettlich2020} as well as further investigated in \cite{Goettlich2021}. On the other hand, accidents are generated due to the self-excitation property of the Hawkes process. These represent for example car pile-ups that may happen due to different reasons. In combination with bad road conditions after rain or snow fall drivers are not able to break their vehicles in time and crash into each other multiple times. This effect does not have to be only due to bad road conditions, it can also be evoked by inattentive or careless drivers right behind an accident site. Another interpretation of the self-excitation property of accidents can be found in traffic jams resulting from a primary accident whose tails may be located on road section that are badly visible like behind a curve. Then again there can be not sufficient space to break in time.  
	The existence of these self-excitation accidents is mentioned in a data analysis in \cite{diLoro.2024,Kalair2020}, stating that 6\% up to 10\% of traffic accidents are of the self-excitation type. We will further quantify the self-excitation in terms of intermediate accident times and finally compare results of a numerical simulation with a set of real data.
	
	This article is now structured as follows. In Section \ref{sec3} we start introducing the traffic accident model for one road and then first discuss how accidents affect traffic flow and second how traffic flow influences the generation of accidents. In the second part, we introduce the self-excitation property, underline it with a data study and define the Hawkes process in the context of accidents. Section \ref{Sec: ExtNetwork} extends the study on one road to a whole traffic network introducing appropriate coupling conditions. The work is concluded in Section \ref{sec: numerics}, where based on a suitable numerical discretization different network examples are studied regarding travel time and routing strategies.  

	\section{Traffic accident model on a single road}\label{sec3}
	In this section we present a traffic accident model that contains self-excitation accidents making use of the Hawkes process. We start introducing the Lighthill-Whitham-Richards (LWR) model \cite{Lighthill1955,Richards1956}, which is a macroscopic traffic model of first order that describes the traffic density on a single road by a hyperbolic conservation law
	\begin{align}
		\label{consLaw}
		\rho_t(x,t) + F(x,t,\rho(x,t))_x = 0, ~~~ \rho(x,0)=\rho_0(x),
	\end{align}
	where $F$ denotes the flux function, $\rho(x,t)$ the traffic density at time $t$ and position $x$ and $\rho_0$ some initial density.
	We extend this model following \cite{Goettlich2020} and include traffic accidents as capacity drops in the flux function. This leads to an additional space dependency in the flux function. For  better illustration, we concentrate on a single road first and extend the framework later to traffic networks. 
	
	\subsection{Influence of accidents on traffic flow dynamics}
	\label{sec2.1}
	We specify more precisely the form of flux function in (\ref{consLaw}) which we decompose into
	\begin{align}
		\label{fluxfunction}
		F(x,t,\rho(x,t)) = c_\text{a}(x,t)c_{\text{road}}(x)f(\rho(x,t)).
	\end{align}
	Typically, traffic flow is characterized by a LWR-type function $f:[0,1] \rightarrow [0,\infty)$, i.e. $f(0)=f(1) = 0$, $f$ being strictly concave. There exists a unique $\rho^* \in (0,1)$ such that $f'(\rho^*)=0$. A common choice is $f(\rho)=\rho v(\rho)=\rho(1-\rho)$, for $v$ the velocity function, see \cite{Garavello2016, Garavello2006}. {We consider a road that is given by the interval $[a,b], ~a<b$ where $a$ and $b$ can take the values $-\infty$ and $\infty$, respectively.} The general road capacity is described by a function $c_{\text{road}}: {[a,b]} \rightarrow \mathbb{R}_{\geq 0}$ and might be for example influenced by the speed limit or the number of available lanes. The function $c_\text{a}: {[a,b]} \times \mathbb{R}_{\geq 0} \rightarrow \mathbb{R}_{\geq 0}$ takes potential accidents into consideration. 
	An accident is described by five parameters: $p$ for the accident position, $s$ for the accident size and $c$ for the accident capacity reduction. Additionally, we equip each accident with an accident time $t$ and a duration $d$ for which the accident lasts. Then, an accident leads to a capacity reduction of $c$ on the road section $[p-\frac{s}{2},p+\frac{s}{2}]$ within the time interval $[t,t + d)$. To capture more than only one accident, the parameters extend to vectors $p \in {[a,b]}^J, s \in \mathbb{R}^J_{\geq 0}, c \in [0,c_{\text{max}}]^J$ for some maximal capacity reduction $0 \leq c_{\text{max}}<1$, $t \in \mathbb{R}^J_{\geq 0}$, $d \in \mathbb{R}^J_{\geq 0}$ and $J \in \mathbb{N}$ accidents. We end up defining the accident capacity function $c_\text{a}$ for $t \in \cap_{j=1}^{J} [t_j,t_j+d_j)$ as
	{\begin{align}
		\label{accidentCapFunction}
		c_\text{a}(x,t) = \prod_{j=1}^J \left(1-c_j \mathbbm{1}_{\left[p_j-\frac{s_j}{2},p_j+\frac{s_j}{2}\right]}(x)\right).\end{align}}
	For better readability, we omit the dependency of $c_\text{a}$ on the accident parameters $p,s$ and $c$. Note that $c_a(x,\cdot)$ is piecewise constant in time for any $x \in {[a,b]}$.
 Therefore, we set up a new initial value problem after a change in $c_a$ happened and therefore consider from now on $c_a$ and the flux function $F$ from  \eqref{fluxfunction} just to be dependent on the space component.
	 We introduce the space of functions with bounded total variation as a suitable solution space.
 \begin{definition}
     A function $\rho: [a,b] \rightarrow \mathbb{R}$ is said to be of bounded total variation if the total variation $TV(\rho,[a,b])$, i.e.
     \begin{align*}
         TV(\rho, [a,b]) = \sup \left\lbrace \sum_{j=1}^N |\rho(x_j) - \rho(x_{j-1})| : a\leq x_0<x_1<\ldots <x_N\leq b, N \in \mathbb{N} \right\rbrace
     \end{align*}
     is finite. The function space containing all functions of bounded total variation on the interval $[a,b]$ is denoted by $BV([a,b])$.
 \end{definition}
Under certain conditions we can ensure the existence of a unique entropy solution to (\ref{consLaw})--(\ref{accidentCapFunction}) as a function in BV$({[a,b]})$. We state these conditions in the following Lemma.
	\begin{lemma}
		\label{LemmaUniqueSol}
		Let $f$ be a LWR-flux, i.e. $f(0)=f(1) = 0$, being strictly concave and assume that there exists a unique $\rho^* \in (0,1)$ such that $f'(\rho^*)=0$. Assume that $c_\text{a}(\cdot,t)c_{\text{road}}(\cdot) \in C^2({[a,b]})\cap TV({[a,b]})$. Furthermore, we assume that $c_\text{a}(\cdot)c_{\text{road}}(\cdot)$, its first derivative, $f$ and $f'$ are essentially bounded and that $c_\text{a}(\cdot)c_{\text{road}}(\cdot)$ and its first derivative are integrable. For the initial condition we require $\rho_0 \in BV({[a,b]})$. 
		Then, there exists a unique entropy solution $\rho(\cdot,t) \in BV({[a,b]})$ to (\ref{consLaw}).
	\end{lemma}
 For a proof we refer to \cite{Goettlich2020}. We remark that in the model~(\ref{consLaw})--(\ref{accidentCapFunction}) there is a dependence of the traffic density on space, time and the density itself. However, in the particular case of the accident framework, we can exploit an additional structure such that the capacity function depends only \emph{piecewise constant} on the temporal component and changes only when an accident occurs or vanishes.
 Therefore, for any of these time intervals, where $c_\text{a}$ is constant, Lemma \ref{LemmaUniqueSol} is applicable and yields the existence of a unique entropy solution. Proceeding in time, we may use the system state at the preceding time horizon as new initial data and consider another problem on the next interval with adjusted capacity function.
 %
 %
 Sequentially, we obtain the existence of solutions on any arbitrary time interval.
	To ensure a solution as described in Lemma \ref{LemmaUniqueSol} we have to remove potential discontinuities in $c_{\text{road}}$ and $c_\text{a}$ from (\ref{accidentCapFunction}) and smooth them out, for example using mollifiers, such that they are at least twice continuously differentiable in space.

	\subsection{Influence of traffic on accidents}
	\label{sec4}
	After having discussed how accidents influence traffic flow, we now present how traffic flow influences accidents in our model. There are various approaches on modeling the accident risk. One could assume a deterministic, constant or time-dependent accident risk or alternatively exploit statistical methods to quantify the accident risk. However, we are interested in a framework that is directly connected to the current traffic situation. Furthermore, we aim to account for the phenomenon of multiple collisions or pile-ups, where several accidents happen in a small time period and a small spatial area. Such accidents resulting out of another will be classified \textit{self-excitation accidents}. 
	
	To validate that the self-excitation property can be found in accident data, we consider intermediate accident times. This is the time that passes between two consecutive accidents in a specified region. It is not easy to detect whether a single accident was caused in direct connection to another accident. Therefore, we do not concentrate on classifying a single accident to be of self-excitation type, but are rather interested in an aggregated view on the intermediate accident times. If accidents happened independently from another with some given rate, we {assume that they are exponentially distributed with some suitable parameter, which might be reasonable from the point of view of queuing theory and can also approximately be extracted from real data, see Figure \ref{intTimes1}.}
	
\begin{figure}[h!]
		\centering
		\includegraphics[scale=0.75]{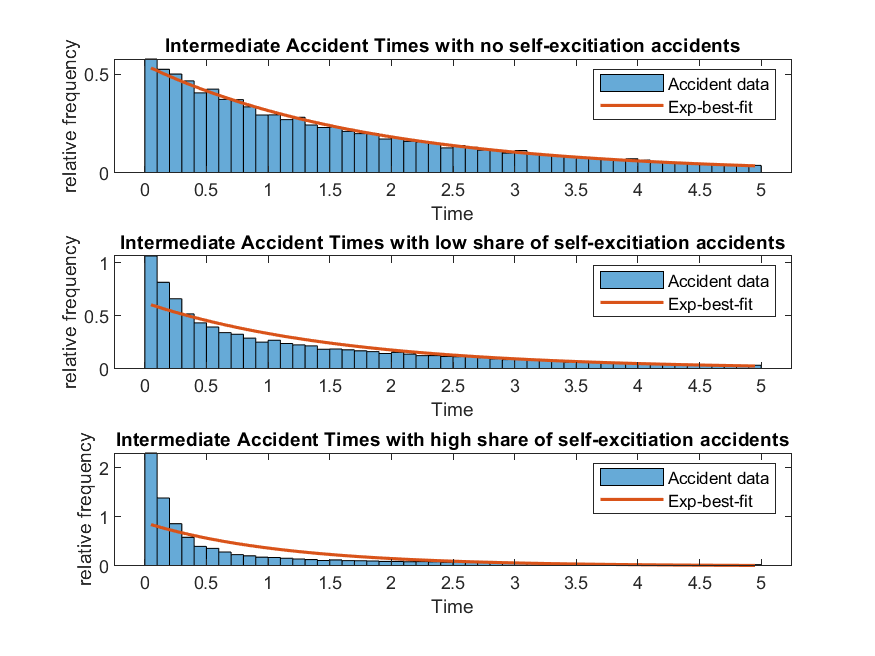}
		\caption{Comparison of intermediate accident times with different shares of self-excitation accidents.}
		\label{intTimes}
	\end{figure}

\begin{figure}[ht!]
		\centering
		\includegraphics[scale=0.59]{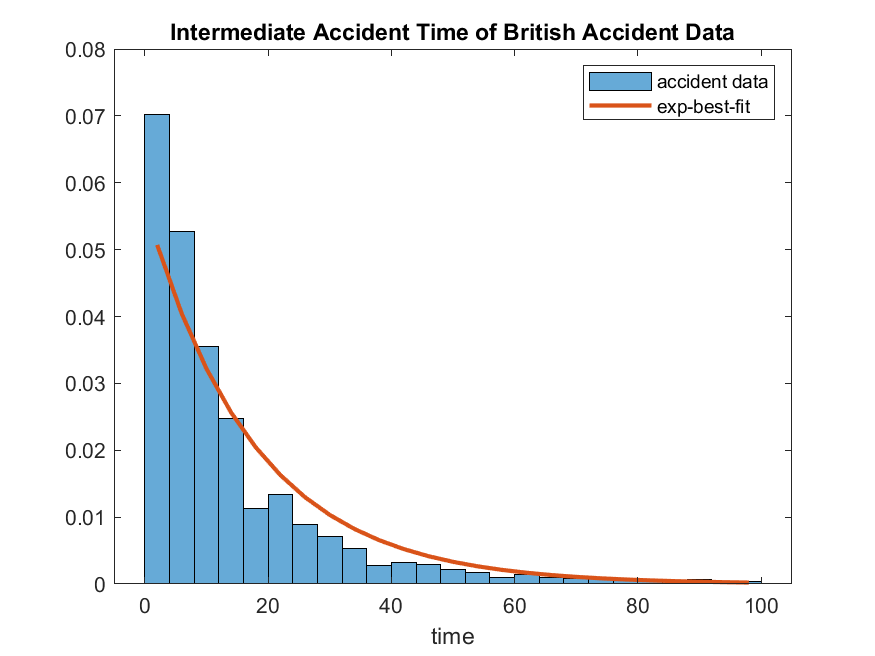}
		\caption{Intermediate accident times for accidents in southern Great Britain in 2018.}
		\label{intTimes1}
	\end{figure}

	In order to underline this property in the following, we consider the diamond network that is studied intensively in Section \ref{sec: numerics} and graphically shown in Figure \ref{networkDiamond}. Its numerical treatment and most parameter choices are set as described in Section \ref{Sec: riskMeasures} and are not of larger interest at this point. We use the diamond network as an artificial traffic network on which we consider three different scenarios, with high 
	and low 
	share of self-excitation accidents and a benchmark setting without a self-excitation component 
	in the accident generation. 
	In Figure \ref{intTimes} we compare the intermediate accident times for the whole diamond network in all three cases.
	
	We can observe that in the benchmark setting the intermediate accident times follow very well the density of an exponential distribution with an appropriate parameter. {When} increasing the self-excitation component, the shares of very low intermediate times increase significantly. As a comparison we again plot the best fit for the density of an exponential distribution and notice that the histogram does not follow the density of an exponential distribution. 
	
	
	To validate which {scale of self-excitation} is the most suitable one, we consider accident data from Great Britain from 2018\footnotemark[1]\footnotetext[1]{\url{https://data.gov.uk/dataset/cb7ae6f0-
			4be6-4935-9277-47e5ce24a11f/road-safety-data},   \\ Accessed: 2024-11-05}. In comparison to many other accident databases, the British database tracks accident times quite exactly and also delivers information on how many vehicles have been involved into an accident. For our purpose, we restrict to accidents lying in a rectangle with a longitude value in $[-1,0]$ and latitude value within $[51,52]$ representing an area in the south-east of London. For this area there were around 30 000 {accidents} reported in 2018. The choice is arbitrary and leads to similar results for different areas with a similar number of accidents. The intermediate accident times are presented in Figure \ref{intTimes1}. We observe that the first three bars (intermediate times up to 6 minutes) exceed all the value for the best exponential fit, whereas times of around 20 minutes are underrepresented. This allows the assumption that there is a self-excitation component in the accident, which is comparable with the artificial generated data from the low share setting in Figure \ref{intTimes}. Summing up so far, this motivates to include the self-excitation effect into the traffic accident model.
	
	Mathematically, we consider a stochastic process which has a self-excitation property to model the times of an accident by its jumps. A suitable stochastic process that has the self-excitation property is called Hawkes process which is introduced in the following and directly related to the accident framework. For a more comprehensive overview on the Hawkes process we refer to \cite{Laub2022}.
	
	The Hawkes process belongs to the class of \textbf{counting processes} {$(N(t),t\geq 0)$ }which are stochastic processes taking values in $\mathbb{N}_0$ and fulfill the following properties:
	\begin{align*}
		{N(0) = 0, ~ N(t)} < \infty \text{ a.s.},~ N(t) \text{ is a pathwise right-continuous step function with stepsize 1}.
	\end{align*}
	{The} sequence of times at which the counting process has a jump $\mathcal{T} = (t_1,t_2,\dots)$, where $\mathcal{T}$ defines a point process which is associated to ${N(t)}$ if $t_1\leq t_2\leq \dots$ holds true almost surely. These jump times represent the times of an accident in our model. Further theoretical investigations for general counting and point processes can be found in \cite{Daley2005}.\\
	To equip the processes with a memory we introduce $(\mathcal{F}_t, t \geq 0)$ the filtration that is generated by $N(t)$ which contains the whole history of information of $N(t)$, especially all previous accidents. Then, for characterizing the jump behaviour, we can define a \textbf{conditional intensity function}. If a non-negative function $\lambda^*(t)$ exists such that 
	\begin{align}
 \label{eq: LimCondExp}
		\lambda^*(t) = \lim_{\Delta t \rightarrow 0}\frac{\mathbb{E}[N(t + \Delta t) - N(t) ~|~ \mathcal{F}_t]}{\Delta t}
	\end{align} then it is called {conditional intensity function}. In the context of traffic accidents $\lambda^*$ represents the likeliness of an accident at time $t$. We are now able to define a Hawkes process.
	\begin{definition}
		\label{def: 2}
		Consider a counting process $(N(t), t\geq 0)$ equipped with the filtration $(\mathcal{F}_t, t\geq 0)$ that satisfies 
		\begin{align*}
			P(N(t+\Delta t) - N(t) = m ~|~ \mathcal{F}_t) = \begin{cases} o(\Delta t), & m>1\\
				\Delta t \lambda^*(t)  + o(\Delta t), & m=1\\
				1 - \Delta t \lambda^*(t)  + o(\Delta t), & m=0.
			\end{cases}
		\end{align*}
		If the conditional intensity function is of the form
		\begin{align*}
			\lambda^*(t) = \lambda(t) + \int_0^t \mu(s)dN(s)
		\end{align*}
		for some $\lambda: [0,\infty) \rightarrow [0,\infty)$ and $\mu:[0,\infty) \rightarrow [0,\infty)$ with $\mu(\cdot)\neq 0$, the process $N(t)$ is called a (linear) \textbf{Hawkes process}. The functions $\lambda$ and $\mu$ are the background intensity and the excitation function, respectively.
	\end{definition}
	The background intensity function models the underlying noise for accidents in the Hawkes process ignoring any clustering effects. Allowing $\lambda$ to be time-dependent, we are able to capture different accident probabilities within a time period and can directly connect them to the traffic situation. For $\gamma>0$, we relate the background accident risk to the flux by
	\begin{align}
		\label{backgroundFct}
		\lambda(\rho(\cdot,t)) = \gamma {\int_{a}^{b}}F(x,\rho(x,t))dx.
	\end{align}
	The choice of $\lambda$ follows two main ideas: The background accident risk on the one hand should increase if the traffic density increases and on the other hand also be positively correlated to the vehicles' velocities. Since $v(\rho) = 1-\rho$ both quantities are negatively correlated and the highest accident risk is located at a certain mixture of high density and high velocity. Therefore, the accident risk is described by the product of both quantities, which is given by the flux function \eqref{fluxfunction}. The choice also leads to the effect that in both extremes, which is $\rho=1$ (maximal density) and $\rho=0$, there is no risk of an accident. In the first scenario, traffic is fully congested and vehicles cannot move. In the second one, there is no traffic at all, which provides a reasonable explanation for the not existing accident risk for both two cases.
	In order to obtain a well-posed Hawkes process, $\lambda$ must remain bounded.
 Under similar assumptions as in Lemma \ref{LemmaUniqueSol}, especially the boundedness of the capacity functions and the flux function, which implies boundedness of the velocity function, we also obtain boundedness of $\lambda$ in time
	\begin{align}
		\begin{split}
			\label{lambdaBeschr}
			|\lambda(\rho(\cdot,t))| &= \left|\gamma\int_{{a}}^{{b}} c_\text{a}(x)c_{\text{road}}(x) \rho(x,t) v(\rho(x,t))dx \right| \\
			&\leq \gamma \|c_\text{a}\|_\infty\|c_{\text{road}}\|_\infty \|v\|_\infty \int_{{a}}^{{b}} \rho(x,t)dx < \infty
		\end{split}
	\end{align}
	
	In contrast, the excitation function models the clustering effects of accidents. To obtain the self-exciting property of the Hawkes process, the excitation function $\mu$ should be decreasing in time (otherwise the process would be called self-regulating). The general impact of the excitation function is proportional to $\int_0^{{t} }\mu(s)ds$. 
The larger $\mu$ is for small arguments, the more likely is an additional jump in the process quite soon after the initial jump, meaning that a second accident is more likely after another accident. If also a significant part of the mass of $\mu$ is located away of the origin, the self-excitation extends to a longer time range.
	In many applications $\mu$ is chosen to be an exponentially decaying function, i.e,
	\begin{align}
		\label{muExp}
		\mu(t) = \alpha e^{-\beta t},
	\end{align}
	for some constants $\alpha, \beta>0$. Increasing $\alpha$ increases the overall risk of an self-excitation accident. The parameter $\beta$ models the shape of temporal decay of the self-excitation accident risk. A larger $\beta$ leads to lower intermediate accident times caused by self-exciting accidents, meaning that the time period between two accidents tends to be smaller. 
	Choosing $\mu$ as in \eqref{muExp}, we can determine a condition that ensures that the Hawkes process stays bounded and we do not end up in a regime with an unlimited number of traffic accidents. Therefore, we introduce the constant
	\begin{align}
		\label{nStar}
		n^* = \int_0^\infty \mu(s)ds = \int_0^\infty \alpha e^{-\beta s} ds = \frac{\alpha}{\beta}.
	\end{align}
	This ratio determines the averaged number of accidents that are caused by one previous accident.\\
	\begin{lemma}
		\label{LemmaNoExplosion}
		Let $\mu$ be given as proposed in (\ref{muExp}) and let $n^*$ from (\ref{nStar}) be strictly smaller than 1. If $\lambda$ is a bounded function for which the limit $\lim\limits_{t \rightarrow \infty} \lambda(t) = \lambda_\infty$ exists, then
		\begin{align*}
			\lim_{t \rightarrow \infty}\mathbb{E}[\lambda^*(t)] = \frac{\lambda_\infty}{1 - n^*} < \infty.
		\end{align*}
	\end{lemma}
	A proof can be found in \cite{Asmussen2003}, Proposition 7.4, where also further details on the interpretation of $n^*$ are provided. {Note that if the conditions in Lemma \ref{LemmaNoExplosion} do not hold, there is a positive probability that the Hawkes process explodes in finite time.}
	If a jump in the Hawkes process on average leads to less than one additional jump, the Hawkes process is well-defined and does not explode if the background noise is also bounded. For our particular setting we therefore assume in the following that $\alpha<\beta$ which is equivalent to $n^*<1$ in (\ref{nStar}).
	
\begin{figure}[htb!]
		\centering
		\includegraphics[scale=0.7]{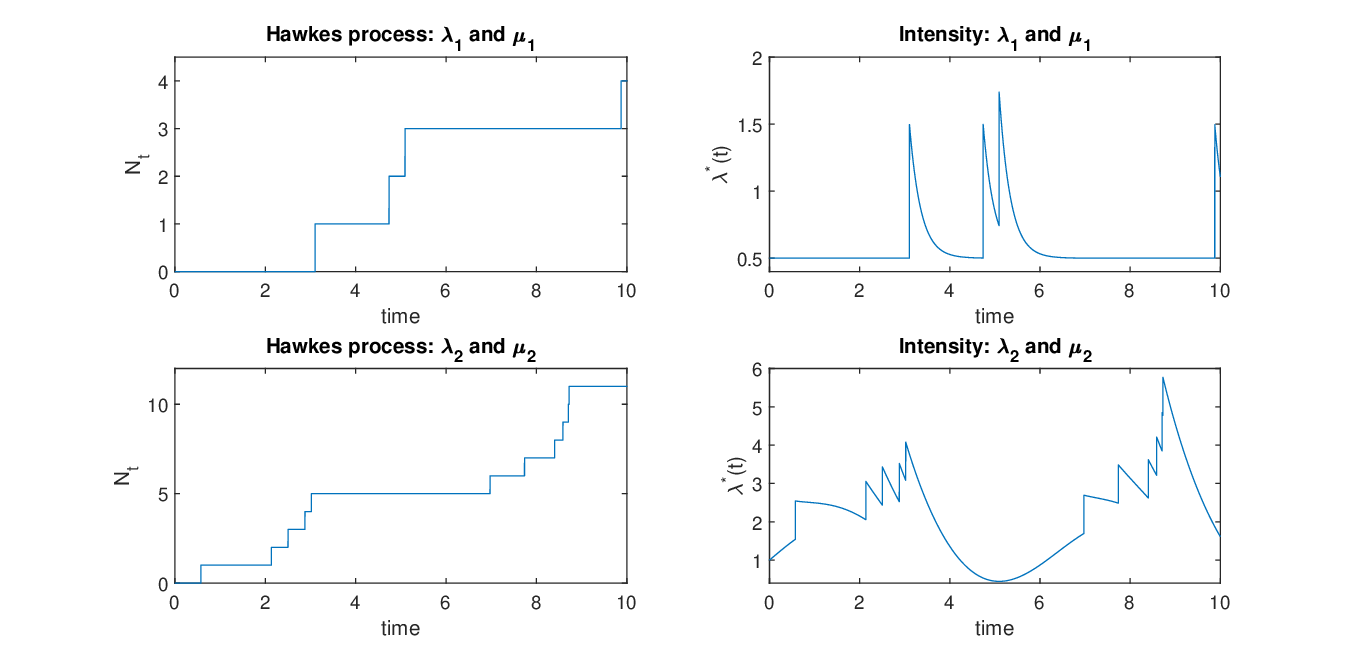}
		\caption{The Hawkes processes $N(t)$ to the left and the corresponding conditional intensity functions $\lambda^*$ to the right, with $\mu_1$, $\lambda_1$ in the first row and $\mu_2$ and $\lambda_2$ in the second row.}
		\label{fig1}
	\end{figure}

	To illustrate an exemplary evolution of the Hawkes process we present in Figure \ref{fig1} two realizations for excitation functions $\mu_1(t) = e^{-4t}, ~ \mu_2(t)=e^{-t}$ and background intensity functions $\lambda_1(t) = \frac{1}{2}, ~\lambda_2(t) = \sin(t) + 1$.
	For $\lambda_1$ and $\mu_1$ we observe that the intensity turns back quite quickly to the base level of $\frac{1}{2}$. Generally, we observe a small amount of jumps and only the third jump seems to be a jump generated by self-excitation. For $\lambda_2$ and $\mu_2$ the background intensity function is sine-shaped and since $\mu_2$ is slower decaying than $\mu_1$, there are more jumps generated than in the first example. Jumps cluster in two areas, around $t=3$ and $t=8$. In between, there is a period without jumps, in which also the conditional intensity function returns to its background sine-shaped part.


	To generate a sequence of jump times of the Hawkes process a so-called thinning procedure is often applied. Such algorithms are typically used to illustrate point processes and have also been applied to the Hawkes process, e.g. in \cite{Ogata1981}. Here, due to the more complex structure of the background intensity function, we proceed differently and evaluate in every time step whether a new accident is generated. We use Definition \ref{def: 2} and a uniformly on $[0,1]$ distributed random variable $u_1 \sim \mathcal{U}([0,1])$ such that a new accident is set if $u_1 \leq \Delta t \lambda^*(t)$. For the algorithm of the entire traffic accident see Algorithm \ref{alg:Hawkes}.

	So far, we only discussed \textit{when} accidents happen in the model. As a next step we describe \textit{where} the accident positions are located under the condition that there is an accident. Therefore, we propose two different approaches, one for the background accidents and one for the accidents generated by self-excitation. 
	%
	For the background accidents we follow \cite{Goettlich2020} using the idea that the likeliness of an accident is positively correlated to both, high traffic density and velocity. The function $\lambda$ defined in (\ref{backgroundFct}) measures this background accident risk.

	For the self-excitation accidents we use another idea. Denote by $P=(p_1, p_2, \dots, p_J)$ the sequence of primary accident positions corresponding to the accident times $\mathcal{T}=(t_1,t_2,\dots, t_J)$. Since we deal with self-excitation accidents the primary accident can only have an influence on the part of the road behind that accident site. Furthermore, we assume that the likeliness for an accident position decreases exponentially in the distance of the primary accident position with some parameter $\tilde{\beta}>0$. Additionally, we assume that there is a small plateau before the exponential decay takes action directly behind the primary accident position. A similar choice has been made in \cite{Kalair2020} based on a data analysis in the UK. The length of the plateau is denoted by $\nu\geq 0$ and according to the data analysis in \cite{Kalair2020} can be estimated to approximately 100 meters. For the primary accident $j$ at $t_j$ with position $p_j$ we define 
	\begin{align}
		\begin{split}
			\label{lambdaTilde}
			\tilde{\lambda}_j(t,B) &= \int_{{a}}^{p_j} \mathbbm{1}_{B}(x) \bigg( \mathbbm{1}_{(p_j - \nu, p_j]}(x) + \mathbbm{1}_{({a}, p_j-\nu]}(x)e^{-\tilde{\beta}(p_j - \nu -x)} \bigg)dx
		\end{split}
	\end{align}
	for $B \in \mathcal{B}({[a,b]})$ a Borel set on ${[a,b]}$. This function describes the likeliness of the accident position for a secondary self-excitation accident for primary accident $j$. We can estimate
	\begin{align}
		\label{boundLambda}
		\left|\tilde{\lambda}_j(t,B) \right| \leq \left|\tilde{\lambda}_j(t,{[a,b]}) \right| \leq  \nu +  \int_0^{b-a} e^{-\tilde{\beta}x}dx = \nu + \frac{1}{\tilde{\beta}} - \frac{e^{-\tilde{\beta}(b-a)}}{\tilde{\beta}} < \infty,
	\end{align}
	which yields the boundedness of $\tilde{\lambda}_j$. We finally define an accident position measure $\mu_t^\text{pos}$ on $({[a,b]}, \mathcal{B}({[a,b]}))$ by
	\begin{align}
 \label{eq: muPosition}
		\mu_t^\text{pos} (B) = \dfrac{\gamma \int_{B} F(x,\rho) dx +  {\frac{\alpha\tilde{\beta}}{\tilde{\beta}\nu + 1 - e^{-\tilde{\beta}(b-a)}}} \sum_{j=1}^J   e^{-\beta(t-t_j)} \tilde{\lambda}_j(t, B) }{\lambda^*(t)}
	\end{align}
	for $B \in \mathcal{B}({[a,b]})$ and $J \in \mathbbm{N}_0$ the number of previously existing accidents. The first term accounts for the background share, whereas the second term contributes a position measure for the self-excitation accidents. The factor $\frac{\alpha\tilde{\beta}}{\tilde{\beta}\nu + 1 - e^{-\tilde{\beta}(b-a)}}$ is a normalizing constant for $\tilde{\lambda}_j$. Additionally, we make use of the fact that 
	\begin{align*}
		\int_0^{{t_J}} \mu(s) dN(s) = \sum_{j=1}^J \mu({t_j}),
	\end{align*}
	where {$N(t)$} has jumps at $t_1,\dots, t_J$.
	
\begin{figure}[ht!]
		\centering
		\begin{minipage}{0.325\textwidth}
			\includegraphics[scale=0.345]{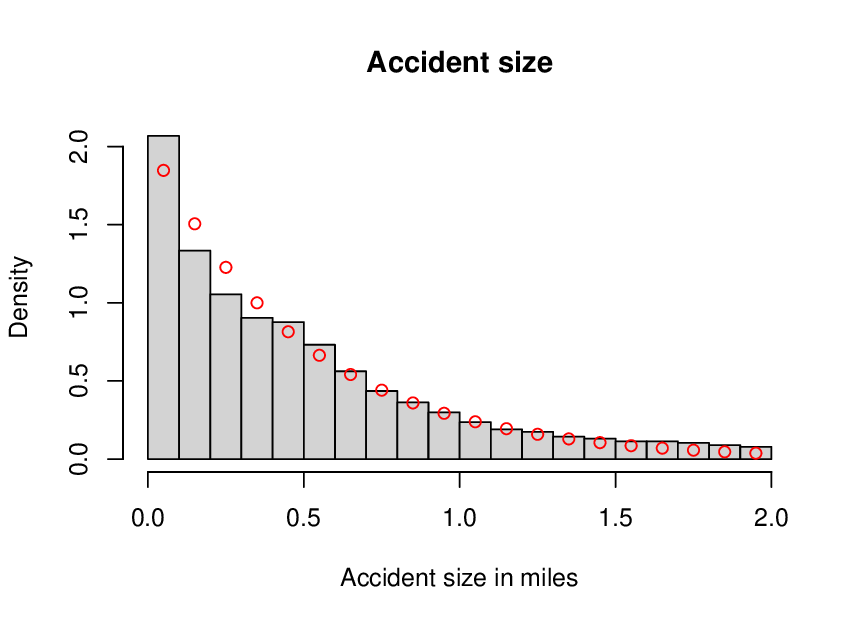}
		\end{minipage}
		\begin{minipage}{0.325\textwidth}
			\includegraphics[scale=0.345]{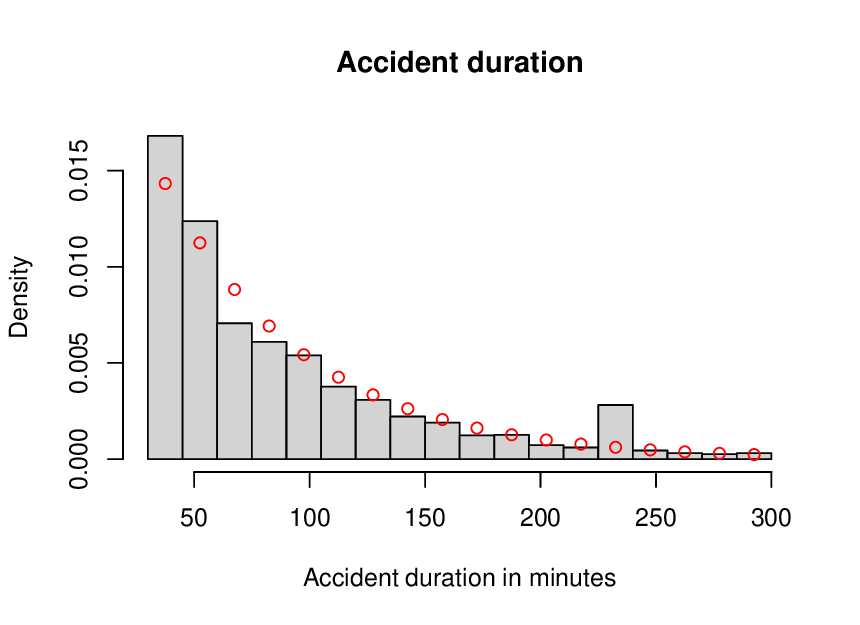}
		\end{minipage}
		\begin{minipage}{0.325\textwidth}
			\includegraphics[scale=0.345]{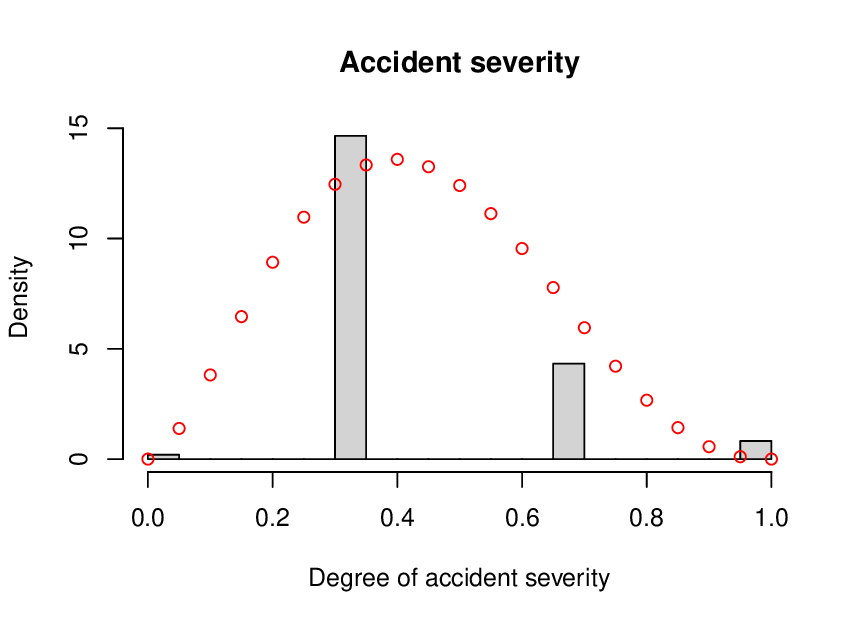}
		\end{minipage}
		\caption{Data comparison of the three parameters compared to an optimal fit for exponential distributions (distance and duration) and Beta distribution (accident severity). Left the accident size, in the middle the accident duration and to the right the accident severity.}
		\label{furtherParams}
	\end{figure}

	It remains to determine probability measures for the accident sizes $s_j$, the accident capacity reductions $c_j$ and the accident duration $d_j$. Even though one could think of a connection of these parameters to particular traffic situations, we treat them independently from the traffic flow, but still as a random variable. The measure for the accident size should provide a positive accident length and it can be assumed that low sizes are more frequent than very long sections where an accident affects the road. Using US accident data of about 3 million accidents since 2016 (see \cite{Moosavi2019}), for a rough approximation, it seems appropriate to choose an exponential distribution for the description of the accident sizes. The parameter for the best fit exponential distribution on the left of Figure \ref{furtherParams} for the US data set can be estimated to be the inverse of the mean accident length of about half a mile or 800 meters.
	
	The capacity reduction is always a real value in $[0,1)$. In real data, we identify the capacity reduction by the accident severity. Considering US accident data from \cite{Moosavi2019} the accident severity is given by four discrete values (from low 0 to high 1). In reality this value might rather be a continuous variable than a discrete one. Therefore, we assume that the capacity reduction is Beta-distributed for some parameters $\hat{\alpha}>0$ and $\hat{\beta}>0$. 
	A parameter estimation on the US data in the right picture of Figure \ref{furtherParams} leads to $\hat{\alpha} = 2.66$ and $\hat{\beta} = 3.53$ which represents a left skewed distribution with mean around 0.43.
	
	It is left to determine the accident duration parameters $d$. We assume that there is some base accident duration $\bar{d}$, which can be explained by the fact that it requires some time until the police arrives and the vehicles are recovered. After $\bar{d}$ elapsed, the additional time the accident is active $d_{\text{add}}$ is assumed to be exponentially distributed with some parameter $\lambda_d$, i.e. $d_{\text{add}} \sim Exp(\lambda_d)$ and $d = \bar{d} + d_{\text{add}}$. These assumptions are supported by the US accident data, where $\bar{d}$ can be estimated to be half an hour and $\lambda_d$ such that the mean of the exponential distribution corresponds to 60 minutes. In the duration plot of Figure \ref{furtherParams} the distribution for $d$ is shown and the origin of the illustration is located in 30 minutes. 
	This concludes the discussion on how traffic affects the generation of accidents in the single road scenario. 

 Summarizing, we present the entire traffic accident model. The dynamics for the traffic density are governed the conservation law in equation \eqref{consLaw}, where the flux functions is chosen as presented in \eqref{fluxfunction}. The traffic accident parameters entering the traffic dynamics via the flux function are then determined as follows: For the accident time, we consider a Hawkes process as in Definition \ref{def: 2} with background intensity $\lambda$ given by equation \eqref{backgroundFct}. The times of an accident $t_j$ are given by the jump times of the Hawkes process. At each of the jump times $t_j$, we determine the four remaining accident parameters $p_j \sim \mu_{t_j}^\text{pos}$ defined in \eqref{eq: muPosition}, $s_j \sim Exp(\lambda_s)$, where $\lambda_s$ is chosen such that the exponential distribution has a mean that corresponds to 800 meters, $c_j \sim Beta(2.62,3.53)$ and $d_j \sim \bar{d} + Exp(\lambda_d)$, where $\bar{d}$ represents 30 minutes and $\lambda_d$ a parameter for the exponential distribution such that its mean is 60 minutes.
 
 Note that fixing one realization for the Hawkes process and the accidents' parameter configurations leads to a piecewise constant accident capacity function $c_\text{a}$ for which we can apply Lemma \ref{LemmaUniqueSol} on any time segment where $c_\text{a}$ is constant in time. Applying Lemma \ref{LemmaUniqueSol} sequentially, we obtain the \textit{pathwise} existence of solutions on an arbitrary time interval.
	
	\section{Traffic accident network model}
	\label{Sec: ExtNetwork}
	As a next step we generalize the framework to a traffic network which we consider to be a directed and finite graph $G = (\mathcal{V},\mathcal{E})$, where $\mathcal{V}$ is a nonempty set of vertices and $\mathcal{E}$ is a nonempty set of arcs modeling the junctions and roads, respectively. Every arc $e \in \mathcal{E}$ is represented by an interval $[a_e,b_e]$ where we set $a_e=-\infty$ for arcs which do not originate in a vertex and $b_e = +\infty$ for arcs which do not discharge into a vertex. On any arc a LWR-type hyperbolic conservation law governs the traffic dynamics presented in (\ref{consLaw}) equipped with the corresponding indices and $e \in \mathcal{E}$ 
	\begin{align}
 \label{eq:TrafficNetwork}
		(\rho_e)_t(x,t) + F_e(x,\rho_e(x,t))_x = 0 ,~~~~ \rho_e(x,0)= (\rho_e)_0(x),
	\end{align}
	where the flux function can be decomposed into
	\begin{align}
 \label{eq:FluxNetwork}
		F_e(x,\rho_e(x,t)) = c_\text{road}^e(x) c_a^e(x) \rho_e(x,t)v_e(\rho_e(x,t)).
	\end{align}
	
	As already pointed out the existence of unique entropy solutions for each single arc $e$ is ensured by Lemma~\ref{LemmaUniqueSol}. However, the network approach additionally requires a discussion on weak solutions at the vertices of the network. To do so, we follow the ideas proposed in \cite{Garavello2006} for the classical LWR model and highlight the key differences.
	
	A \textit{weak solution} at a vertex $v \in \mathcal{V}$ with $n$ ingoing and $m$ outgoing roads represented by intervals $I_i^{\text{in}} = [a_i^{\text{in}},b_i^{\text{in}}]$ and $I_j^{\text{out}} = [a_j^{\text{out}},b_j^{\text{out}}]$ is a collection of functions $\rho_i^{\text{in}}: I_i^{\text{in}} \times \mathbb{R}_{\geq 0} \rightarrow \mathbb{R}, ~ i=1,\dots n$ and $\rho_j^{\text{out}}: I_i^{\text{out}} \times \mathbb{R}_{\geq 0} \rightarrow \mathbb{R}, ~ j=1,\dots m$ such that
	\begin{align}
		\begin{split}
			\label{weakSolution}
			&\sum_{i=1}^n \left(\int_0^\infty \int_{I_i^{\text{in}}} \rho_i^\text{in}(x,t) \partial_t \phi_i^\text{in}(x,t)  + F_i^{\text{in}}(x,\rho_i^\text{in}(x,t)) \partial_x \phi_i^\text{in}(x,t) dx dt  \right)\\
			+~& \sum_{j=1}^m \left(\int_0^\infty \int_{I_j^{\text{out}}} \rho_j^\text{out}(x,t) \partial_t \phi_j^\text{out}(x,t)  + F_j^{\text{out}}(x,\rho_j^\text{out}(x,t)) \partial_x \phi_j^\text{out}(x,t) dx dt  \right) = 0,
		\end{split}
	\end{align}
	for every choice of test functions $\phi_i^{\text{in}} \in C_0^1(I_i^{\text{in}} \times \mathbb{R}_{\geq 0}),~ i=1,\dots n$ and $\phi_j^\text{out} \in C_0^1(I_j^{\text{out}} \times \mathbb{R}_{\geq 0}),~ j=1,\dots m$ which satisfy
	\begin{align*}
		\phi_i^\text{in}(b_i^{\text{in}},t) = \phi_j^\text{out}(a_j^{\text{out}},t),~~~~\partial_x \phi_i^{\text{in}}(b_i^{\text{in}},t) = \partial_x \phi_j^\text{out}(a_j^{\text{out}},t),
	\end{align*}
	for all $t\geq 0$. From the definition of the weak solution one can directly deduce that such a collection of functions allows for flux conservation at the vertex $v \in \mathcal{V}$, i.e.,
	\begin{align}
 \label{eq:FluxConservation}
		&\sum_{i=1}^n  F_i^{\text{in}}(b_i^{\text{in}}+,\rho_i^\text{in}(b_i^{\text{in}}+,t)) = \sum_{j=1}^m  F_j^{\text{out}}(a_j^{\text{out}}-,\rho_j^\text{out}(a_j^{\text{out}}-,t)).
	\end{align}
	The expressions $b_i^\text{in}+$ and $a_j^\text{out}-$ denote the limits from the left and right, respectively {and $F_i^{\text{in}}(b_i^{\text{in}}+,\rho_i^\text{in}(b_i^{\text{in}}+,t)),F_j^{\text{out}}(a_j^{\text{out}}-,\rho_j^\text{out}(a_j^{\text{out}}-,t))$ the corresponding in- and outgoing fluxes at the junction}. For the consideration of well-posedness in the following we always assume that $F_i^\text{in}(\lim_{x \nearrow b_i^\text{in}} x, \rho)$ and $F_j^\text{out}(\lim_{x \searrow a_j^\text{out}} x, \rho)$ exist requiring that $F$ has no spatial discontinuity arbitrarily close to the junction. Therefore, we have to state a technical assumption as follows: For all accidents $(p,s,c,t,d)$ there exists $\varepsilon>0$ such that for every vertex $v \in \mathcal{V}$ and the corresponding in- and outgoing roads it holds
	\begin{align}
		\begin{split}
			\label{conditionDiscontinuities}
			&p - \frac{s}{2} \notin \left(\bigcup_{i=1}^n \left[b_i^\text{in} - \varepsilon, b_i^\text{in}\right] \right) \bigcup \left(\bigcup_{j=1}^m \left[a_j^\text{out}, a_j^\text{out} + \varepsilon\right]\right),\\
			&p + \frac{s}{2} \notin \left(\bigcup_{i=1}^n \left[b_i^\text{in} - \varepsilon, b_i^\text{in}\right]\right) \bigcup \left(\bigcup_{j=1}^m [a_j^\text{out}, a_j^\text{out} + \varepsilon]\right).
		\end{split}
	\end{align}
	The technical condition \eqref{conditionDiscontinuities} excludes jumps in the flux function in the spatial component arbitrarily close to the junctions. Conversely, we can also find a (small) road segment with positive Lebesgue measure leading into (or out of) the junction on which the flux function is spatially constant such that results on the existence of solution at the junction~\cite{Garavello2016,Garavello2006} are applicable.
 
	Flux conservation does not necessarily describe the proportion of flows from ingoing to outgoing roads. Thus, we introduce a distribution matrix $A\in \mathbb{R}^{n \times m}$ where $A_{ij}$ is the share of traffic flow leading into the junction from the ingoing road $i$ and flowing into the outgoing road $j$. Flux conservation is ensured if it holds $\sum_{j = 1}^m A_{ij} = 1$ for all $i=1,\dots,n$. The outgoing flows are then given by
	\begin{align}
		\label{admissible3}
		F_j^{\text{out}}(a_j^{\text{out}}-,\rho_j^\text{out}(a_j^{\text{out}}-,t)) = \sum_{i=1}^n A_{ij} F_i^{\text{in}}(b_i^{\text{in}}+,\rho_i^\text{in}(b_i^{\text{in}}+,t)).
	\end{align}
	
	Collecting these conditions and assuming $\rho_i^{\text{in}}(\cdot, t) \in BV([a_i^{\text{in}},b_i^{\text{in}}]), \rho_j^{\text{out}}(\cdot, t) \in BV([a_j^{\text{in}},b_j^{\text{in}}])$ for all $i=1,\dots,n$, $j=1,\dots,m$ and {for all $ t\geq 0$},  
	we can define an \textit{admissible solution} at a vertex $v \in \mathcal{V}$ for a collection of functions $\rho$ if
	\begin{enumerate}
		\item $\rho$ is a weak solution according to $(\ref{weakSolution})$
		\item the condition ($\ref{conditionDiscontinuities}$) is assumed to hold
		\item the fluxes are arranged according to the distribution matrix as in ($\ref{admissible3}$)
		\item $\sum_{i=1}^n F_i^{\text{in}}(b_i^{\text{in}}+,\rho_i^\text{in}(b_i^{\text{in}}+,t))$ is maximal subject to the item 3.
	\end{enumerate}
 {The fourth condition is required to ensure the well-posedness and uniqueness of solutions at the junctions.
        We refer to the books of Garavello and Piccoli~\cite{Garavello2016,Garavello2006} for further details.
        The intuition behind flux maximization with respect to conditions 1.-3. is that any driver crosses the junction if there is space. Therefore, as many vehicles as possible cross the junction. If condition 4 would not be enforced, vehicles may stop in front of an empty junction, which is not reasonable.}
	In the following, we restrict to special networks described by 1-2 ($n=1, m=2$) and 2-1 ($n=2, m=1$) junctions and discuss the latter maximization problem in detail in order to find an admissible solution at the vertex.\\
	\paragraph{1-2 junction}
	\label{sec 1-2junction}
	We employ the demand and supply formulation introduced by Lebacque \cite{Lebacque2005} to determine the maximal fluxes at every vertex and thereby  ensure flux conservation. Here, demands representing the amount of flux that the ingoing road demands to provide, and supplies given by the amount of flux that the outgoing roads are able to receive, i.e.,
	\begin{align*}
		D_i(\rho) = \begin{cases}
			F^{\text{in}}_i(b_i^{\text{in}},\rho), & \text{if } \rho \in [0,\rho^*]\\
			F^{\text{in}}_i(b_i^{\text{in}},\rho^*), & \text{if } \rho \in (\rho^*,1]
		\end{cases}, ~~~ 
		S_j(\rho) = \begin{cases}
			F^{\text{out}}_j(a_j^{\text{out}},\rho^*), & \text{if } \rho \in [0,\rho^*]\\
			F^{\text{out}}_j(a_j^{\text{out}},\rho), & \text{if } \rho \in (\rho^*,1]
		\end{cases},
	\end{align*}
	where $\rho^* \in (0,1)$ is the unique maximizing density of the flux function. Then, the maximization problem motivated by the conditions of an admissible solution reads
	\begin{align}
		\begin{split}
			\label{opt12}
			& \max_{F_1^{\text{in}}(b_1^\text{in}, \rho_1^\text{in})} F_1^{\text{in}}(b_1^\text{in}, \rho_1^\text{in}) \\
			\text{s.t.} ~ &F_j^{\text{out}}(a_j^\text{out}, \rho_j^\text{out}) \in [0,S(\rho_j^\text{out})],~~ j=1,2\\
			&F_1^{\text{in}}(b_1^\text{in}, \rho_1^\text{in}) \in [0,D(\rho_1^\text{in})]\\
			& F_j^\text{out}(a_j^\text{out}, \rho_j^\text{out}) = A_{1j} F_1^\text{in}(b_1^\text{in}, \rho_1^\text{in}),~~ j=1,2
		\end{split}
	\end{align}
	The solution to this optimization problem in the 1-2 network setting is given by
	\begin{align*}
		F_1^{\text{in}}(b_1^\text{in}, \rho_1^\text{in}) = \min \left\lbrace D_1(\rho_1^\text{in}), \frac{S_1(\rho_1^\text{out})}{A_{11}}, \frac{S_2(\rho_2^\text{out})}{A_{12}}  \right\rbrace.
	\end{align*}
	So far only the fluxes have been considered at the junction. It is left to fix the corresponding densities there. However, due to the concave structure of the flux function, there exist typically two different density values to a given flux value.
	Therefore, we follow the literature in \cite{Garavello2016, Garavello2006} and emphasize again that condition ($\ref{conditionDiscontinuities}$) holds. Let us define the mapping $\tau: [0,1] \rightarrow [0,1]$ such that
	\begin{align*}
		f(\tau(\rho)) = f(\rho), ~~ \forall \rho \in [0,1]~~ \text{and} ~~ \tau(\rho) \neq \rho, ~~ \forall \rho \in [0,1] \backslash \lbrace \rho^* \rbrace.
	\end{align*}
	The densities at the junction are then defined by
	\begin{align}
		\begin{split}
			\label{tauRho}
			\rho(b_i^\text{in},t) &= \begin{cases}
				\lbrace\rho_i^\text{in} \rbrace \cup (\tau(\rho_i^\text{in}),1], & 0 \leq \rho_i^\text{in} \leq \rho^*\\
				[\rho^*, 1], & \rho^* < \rho_i^\text{in} \leq 1
			\end{cases}\\
			\rho(a_j^\text{out},t) &= \begin{cases}
				[\rho^*, 1] , & 0 \leq \rho_j^\text{out} \leq \rho^*\\
				\lbrace \rho_j^\text{out}\rbrace \cup [0,\tau(\rho_j^\text{out})), & \rho^* < \rho_j^\text{out} \leq 1.
			\end{cases}
		\end{split}
	\end{align}

	\paragraph{2-1 junction}
	\label{sec 2-1junction}
	Next, we consider a network with two ingoing roads and one outgoing road. Note that in this case, there are no distribution parameters $A_{ij}$ necessary. The optimization problem for the flux maximization reads
	\begin{align}
		\begin{split}
			\label{opt21}
			& \max F_1^{\text{in}}(\rho(b_1^\text{in}),b_1^\text{in}) +  F_2^{\text{in}}(\rho(b_2^\text{in}),b_2^\text{in}) \\
			\text{s.t.} ~ &F_1^{\text{out}}(\rho_1^\text{out},a_1^\text{out}) \in [0,S(\rho_1^\text{out})]\\
			&F_i^{\text{in}}(\rho_i^\text{in},b_i^\text{in}) \in [0,D(\rho_i^\text{in})],~~i=1,2\\
			& F_1^\text{out}(\rho_1^\text{out},a_1^\text{out}) =  F_1^\text{in}(\rho_1^\text{in},b_1^\text{in}) + F_2^\text{in}(\rho_2^\text{in},b_2^\text{in}).
		\end{split}
	\end{align}
	The solution to ($\ref{opt21}$) is not unique if the sum of demands exceeds the supply. Therefore, a so-called rightway parameter $q \in [0,1]$ is introduced that models the priorities of the two ingoing roads. If the demands exceed the supply, i.e. $D(\rho_1^\text{in}) + D(\rho_2^\text{in}))>S(\rho_1^\text{out})$, there are basically three scenarios in which the inflows are set as follows:
	\begin{align*}
		F_1^\text{in}(\rho_1^\text{in},b_1^\text{in}) = 
		\begin{cases}
			qS(\rho_1^\text{out}) & D(\rho_1^\text{in})>qS(\rho_1^\text{out}),~D(\rho_2^\text{in})>(1-q)S(\rho_1^\text{out})\\
			S(\rho_1^\text{out}) - D(\rho_2^\text{in}) & D(\rho_1^\text{in})>qS(\rho_1^\text{out}),~D(\rho_2^\text{in})\leq(1-q)S(\rho_1^\text{out})\\
			D(\rho_1^\text{in}) & D(\rho_1^\text{in})\leq qS(\rho_1^\text{out}),~D(\rho_2^\text{in})>(1-q)S(\rho_1^\text{out})
		\end{cases}
	\end{align*}
	\begin{align*}
		F_2^\text{in}(\rho_2^\text{in},b_2^\text{in}) = 
		\begin{cases}
			(1-q)S(\rho_2^\text{out}) & D(\rho_1^\text{in})>qS(\rho_1^\text{out}),~D(\rho_2^\text{in})>(1-q)S(\rho_1^\text{out})\\
			D(\rho_2^\text{in}) & D(\rho_1^\text{in})>qS(\rho_1^\text{out}),~D(\rho_2^\text{in})\leq(1-q)S(\rho_1^\text{out})\\
			S(\rho_1^\text{out}) - D(\rho_1^\text{in}) & D(\rho_1^\text{in})\leq qS(\rho_1^\text{out}),~D(\rho_2^\text{in})>(1-q)S(\rho_1^\text{out}).
		\end{cases}
	\end{align*}
	If the sum of the demands (the total inflow) does not exceed the supply, the flow at the junction is just given by the demand again. Traffic densities at the junction can be obtained as already presented for the 1-2 junction setting in equation (\ref{tauRho}).
	
	These two particular examples represent an explicit way to define densities and flow at a fixed vertex. Combined with the assumptions of Lemma \ref{LemmaUniqueSol} for every road we obtain a unique solution at a junction $v \in \mathcal{V}$. 

{
 \begin{remark}
     General traffic networks might consist of n-1, 1-m or n-m junctions ($m,n \in \mathbbm{N}$). However, the two prototypical examples for 1-2 and 2-1 junctions show the main ideas and allow for an explicit representation of solutions at the intersection. {For general types of junctions the solution is typically not given explicitly and thus a linear program with the forth condition as objective function and conditions 1.-3. as constraints needs to be solved. As an example, additional theoretical constraints on the rightway parameter are required when ensuring the unique solvability at n-1 junctions.} For more details we refer to \cite{Garavello2016,Garavello2006}.
 \end{remark}
	}
	
	\subsection{Influence of accidents on traffic flow on networks}
	We extend the approaches from Section \ref{sec2.1} to the network setting. Since the roads may have finite length, it is possible that an accident affects more than one road only. We assume that accidents exceeding one road lap over to all in- or outgoing roads. For an accident on road $e \in \mathcal{E}$ we therefore have to adjust the accident capacity functions $c_a^e$ for road $e$ and the in- and outgoing roads of $e$. Denote by $E^\text{in} \subseteq \mathcal{E}$ the set of all ingoing roads into road $e$ and by $E^\text{out} \subseteq \mathcal{E}$ the set of all outgoing roads from $e$. Then, we set for any $e^\text{in} \in E^\text{in}$, $e^\text{out} \in E^\text{out}$ and an accident $(p_j,s_j,c_j,t_j,d_j)$ that is located on road $e$ at $p_j \in (a_e,b_e)$
	\begin{align}
		\begin{split}
			\label{networkCa}
			c_a^{e^\text{in}}(x) &= \begin{cases}
				1,& a_e - p_j + \frac{s_j}{2}\leq 0\\
				1 - c_j \mathbbm{1}_{\left[\max \left\lbrace a_{e^\text{in}}, b_{e^\text{in}} - \left(a_e - p_j + \frac{s_j}{2}\right) \right\rbrace,b_{e^\text{in}}\right]}(x),& a_e - p_j + \frac{s_j}{2}>0
			\end{cases} \\
			c_a^{e}(x) &= 1 - c_j\mathbbm{1}_{\left[\max\left\lbrace a_e, p_j - \frac{s_j}{2}\right\rbrace, \min\left\lbrace b_e, p_j + \frac{s_j}{2}\right\rbrace\right]}(x)\\
			c_a^{e^\text{out}}(x) &= \begin{cases}
				1 , & b_e - p_j - \frac{s_j}{2}\geq 0\\
				1 - c_j \mathbbm{1}_{\left[a_{e^\text{out}}, \min \left\lbrace b_{e^\text{out}}, a_{e^\text{out}} + \left(p_j + \frac{s_j}{2} - b_e\right)\right \rbrace \right]}(x) , & b_e - p_j - \frac{s_j}{2}< 0.
			\end{cases}
		\end{split}
	\end{align}
	If additionally for $e^\text{in}$ it holds $a_{e^\text{in}} - ( b_{e^\text{in}} - \left(a_e - p_j + \frac{s_j}{2}\right))>0$, then the accident on road $e$ also has an influence on the ingoing roads of $e^\text{in}$. This procedure can iteratively be continued until the accident size is reached. The same holds for an accident that exceeds an outgoing road $e^\text{out}$. For more than one accident, the accident capacity function is built by multiplying the accident effects as in the single road case, compare ($\ref{accidentCapFunction}$).

	\subsection{Influence of traffic on accidents on networks}
	This section expands Section \ref{sec4} to networks. First, we consider the choice of the accident times for a traffic network which will be again modeled using a Hawkes process denoted by $(N(t))_{t \in [0,T]}$ and counts the number of accidents that happened starting at $t=0$. Each element of the corresponding point process $\mathcal{T}$ of jump times $\mathcal{T}=(t_1, t_2,\dots)$ describes the time of one accident. Note that there is one Hawkes process for the whole network and not one process for each road. The conditional intensity function is then given by
	\begin{align}
		\label{eq:condIntensity2Network}
		\lambda^*(t) = \left(\sum_{e \in \mathcal{E}}  \lambda_e(\rho_e(\cdot,t)) \right) + \int_0^t \mu(t-s)dN(s).
	\end{align}
	The background intensity functions $\lambda_e$ are defined similar as in (\ref{backgroundFct}) but accounting for a possibly finite length by 
	\begin{align*}
		\lambda_e(\rho_e(\cdot,t)) = \gamma \int_{a_e}^{b_e}F_e(x,\rho_e(x,t))dx.
	\end{align*}
	The functions $\lambda_e$ shall be bounded and the self-excitation function $\mu$ is chosen exactly as in Section \ref{sec4}.
	
	For the accident positions we also have to make some adjustments, since the risk of self-excitation accidents has to be transferred to roads behind the accident. Denote $|E|$ the total number of roads in the traffic network $\mathcal{E}$ and let $\mathcal{R}_j$ for an accident $j$ be the set of roads that are located behind the road on which accident $j$ happens. Additionally, set for $r \in \mathcal{R}_j$, $\kappa_{r,j}$ to the distance between the end of road $r$ and the beginning of the road on which accident $j$ happens. 
	In the network case, we redefine the functions $\tilde{\lambda}_j$ from (\ref{lambdaTilde}) by
	\begin{align*}
		\begin{split}
			\tilde{\lambda}_j(t,B_1,\dots, B_{|E|}) &= \int_{a_r}^{p_j} \mathbbm{1}_{B_r}(x) \bigg( \mathbbm{1}_{[\max\lbrace p_j - \nu,a_r \rbrace, p_j]}(x) \\
			+& \mathbbm{1}_{(a_r,\max \lbrace p_j-\nu,a_r \rbrace]}(x)e^{-\tilde{\beta}(\max\lbrace p_j - \nu,a_r \rbrace -x)} \bigg)dx \\
			+& \sum_{s \in \mathcal{R}_j} \zeta_{s,j}\int_{a_s}^{b_s} \mathbbm{1}_{B_s}(x) \bigg(\mathbbm{1}_{[\max \lbrace b_s - \max \lbrace 0, a_r - (p_j - \nu) - \kappa_{s,j} \rbrace,a_s \rbrace,b_s]}(x)\\
			+& \mathbbm{1}_{[a_s,\max \lbrace b_s - \max \lbrace 0, a_r - (p_j - \nu) - \kappa_{s,j} \rbrace,a_s \rbrace]}(x) e^{-\tilde{\beta}((\max\lbrace p_j - \nu,a_r \rbrace-a_r) + \kappa_{s,j} + (b_s-x))} \bigg)dx
		\end{split}
	\end{align*}
	for $B_e \in \mathcal{B}([a_e,b_e])$ a Borel set over the interval $[a_e,b_e]$.  This more general definition allows for accidents affecting more than one road. The parameter $\zeta_{s,j}$ is a branching factor that decreases the likeliness of an accident on a road in $\mathcal{R}_j$ if there is more than one ingoing road. It is chosen such that $\tilde{\lambda}_j(t,[a_1,b_1],\dots,[a_{|E|},b_{|E|}]) = \nu + \frac{1}{\tilde{\beta}}$. Note that despite the choice of closed intervals $[a_e,b_e]$, accidents almost surely occur inside the road and not at the boundary since $\lbrace a_e \rbrace$ and $\lbrace b_e \rbrace$ are null sets with respect to the Lebesgue integration. 
 The treatment of accidents at the junction is provided in Section \ref{sec: accidentsAtJunction}. 
	
	To determine an accident position in case that there is an accident, the procedure is as follows. First, we select the road of the network on which we generate the accident and then in a second step the exact position $p_j$ is determined.  
	
	We define a probability measure on $(\lbrace 1,\dots,|E| \rbrace, \mathcal{P}(\lbrace 1,\dots,|E| \rbrace))$ selecting the road on which the accident happens by
	\begin{align}
		\label{accmeasChoice}
		\mu_t^\text{RI}(\mathcal{I}) = \dfrac{\left(\sum_{i \in \mathcal{I}} \lambda_{i}(t) \right) + \sum_{j = 1}^J \sum_{i \in \mathcal{I}} \tilde{\lambda}_j(t, \emptyset,\dots,\emptyset, [a_{i},b_{i}],\emptyset,\dots,\emptyset) \alpha e^{-\beta(t-t_j)}}{\lambda^*(t)},
	\end{align}
	for $\mathcal{I} \in \mathcal{P}(\lbrace 1,\dots,|E| \rbrace)$ the power set of $\lbrace 1,\dots,|E| \rbrace$ and where we sum over all primary accidents $j=1,\dots,J$. The boundedness of the background functions has been shown in (\ref{lambdaBeschr}) and the boundedness of $\tilde{\lambda}$ choosing $\zeta_{s,j}$ appropriately is a direct consequence of (\ref{boundLambda}). For a finite number of primary accidents the measure $\mu_t^\text{RI}$ is well-defined. 
	
	It remains to provide an accident position measure on the particularly chosen road. Assume that the accident happens on some road $e \in \mathcal{E}$, then we define a measure $\mu_t^\text{pos,e}$ on $([a_e,b_e], \mathcal{B}([a_e,b_e]))$ by
	\begin{align}
 \label{eq: mute}
		\mu_t^\text{pos,e} (B_e) = \dfrac{\gamma \int_{B_e} F_e(x,\rho_e) dx + \alpha \frac{\tilde{\beta}\nu + 1}{\tilde{\beta}}\sum_{j=1}^J \tilde{\lambda}_j(t, \emptyset,\dots,\emptyset, B_e,\emptyset,\dots,\emptyset) e^{-\beta(t-t_j)}} {\lambda_e(\rho_e(x,t)) + \alpha \frac{\tilde{\beta}\nu + 1}{\tilde{\beta}}\sum_{j=1}^J \tilde{\lambda}_j(t, \emptyset,\dots,\emptyset, [a_e,b_e],\emptyset,\dots,\emptyset) e^{-\beta(t-t_j)}},
	\end{align}
	for $B_e \in \mathcal{B}([a_e,b_e]))$.

	\subsection{Accidents at the junction}
	\label{sec: accidentsAtJunction}
	From our own experience we know that junctions are a location, where accidents are very likely to happen. According to the German Unfallatlas 2020\footnotemark[1] \footnotetext[1]{\url{https://unfallatlas.statistikportal.de/_opendata2022.html, }{Accessed: 2023-11-05}}about 27 percent of the accidents happened with inturning or crossing vehicles. In the US data set from \cite{Moosavi2019} the share of accidents at a crossing or a junction amounts to about 17 percent. Therefore, we extend the introduced traffic accident network model and add an additional component that models accidents directly at the junction. Note that in the case an accident is observed at a junction, the accident capacity functions of the in- and outgoing roads are all affected. For any ingoing road $e^\text{in} \in E^\text{in}$ given by $[a_{e^\text{in}}, b_{e^\text{in}}]$, any outgoing road $e^\text{out} \in E^\text{out}$ given by $[a_{e^\text{out}}, b_{e^\text{out}}]$ and an accident $(p_j,s_j,c_j,t_j,d_j),$ we have
	\begin{align*}
		\begin{split}
			c_a^{e^\text{in}}(x) = 1 - c_j \mathbbm{1}_{\left[\max \left\lbrace a_{e^\text{in}}, b_{e^\text{in}} - \frac{s_j}{2} \right\rbrace,b_{e^\text{in}}\right]}(x),\\
			c_a^{e^\text{out}}(x)=  1 - c_j \mathbbm{1}_{\left[a_{e^\text{out}}, \min \left\lbrace b_{e^\text{out}}, a_{e^\text{in}} + \frac{s_j}{2}\right\rbrace \right]}(x).
		\end{split}
	\end{align*}
	If the accidents overlap the direct in- or outgoing roads we have to further adapt the accident capacity functions of their in- or outgoing functions as presented in (\ref{networkCa}). It should also be remarked that this construction does not contradict the assumption \eqref{conditionDiscontinuities} 
	since the discontinuities in the capacity function occur at $p_j \pm \frac{s_j}{2}$ and not at the accident positions itself. 
	
	It remains to determine the likeliness of such an accident at the junction. Like the assumption on a road, we assume that the accident risk is directly related to the flux at the junction. 
	For a junction $v \in \mathcal{V}$ we denote the flux by $F_v(\rho_1^\text{in}(t),\cdots,\rho_{n_v}^\text{in}(t),\rho_1^\text{out}(t),\cdots,\rho_{m_v}^\text{out}(t))$ and adjust the background accident risk of the Hawkes process for a constant $\gamma_v>0$ to
	\begin{align}
		\label{backgroundFctJunctions}
		\lambda(\rho(\cdot,t)) &= \gamma \sum_{e \in \mathcal{E}}\int_{a_e}^{b_e}F_e(x,\rho_e(x,t))dx + \sum_{v \in \mathcal{V}} \gamma_v F_v(\rho_1^\text{in}(t),\dots,\rho_{n_v}^\text{in}(t),\rho_1^\text{out}(t),\dots,\rho_{m_v}^\text{out}(t)).
	\end{align}
	Then again, the jumps of the Hawkes process with background risk function as given in \eqref{backgroundFctJunctions} models the timings of the accidents.
	After having determined the accident time, we also have to adjust the position probability measures for the network. In particular, we extend $\mu_t^{RI}$ from (\ref{accmeasChoice}), where the indices for the road of the accident are chosen by a component that allows to also choose a junction. The measure acts on $(\lbrace1,\dots,|E|, |E|+1, \dots, |E| + |\mathcal{V}| \rbrace, \mathcal{P}(\lbrace1,\dots,|E|, |E|+1, \dots, |E| + |\mathcal{V}| \rbrace)) $, where the indices $1,\dots, |E|$ represent a road and $|E|+1,\dots,|E|+|\mathcal{V}|$ the junctions, i.e.,
	\begin{align}
		\label{accmeasChoiceJunction}\begin{split}
		     \mu_t^\text{RI}(\mathcal{I}) &= \frac{1}{\lambda^*(t)} \bigg(\gamma \sum_{e \in \mathcal{I}\cap\mathcal{E}}\int_{a_e}^{b_e}F_e(x,\rho_e(x,t))dx\\
  &~~~+ \sum_{v \in \mathcal{I}\cap\mathcal{V}} \gamma_v F_v(\rho_1^\text{in}(t),\dots,\rho_{n_v}^\text{in}(t),\rho_1^\text{out}(t),\dots,\rho_{m_v}^\text{out}(t))   \bigg).
		\end{split}
	\end{align}
 Concluding, the accident network model consists of traffic dynamics on each road given by the conservation law in \eqref{eq:TrafficNetwork} with flux function as in \eqref{eq:FluxNetwork}. At the junctions, we require flux conservation as given in \eqref{eq:FluxConservation}. For the time points of accidents we consider one Hawkes process for the entire network with conditional intensity function $\lambda^*(t)$ presented in \eqref{eq:condIntensity2Network}, where the background intensity in the most general model is given by \eqref{backgroundFctJunctions}. For the accident position, we employ \eqref{accmeasChoiceJunction} for the road or junction index and \eqref{eq: mute} for the precise position on the road. The remaining accident parameters are chosen as in the single lane case. Finally, the accident capacity function with potentially road overlapping accidents is given by \eqref{networkCa} taking into account condition \eqref{conditionDiscontinuities}.

	\section{Numerical experiments}
	\label{sec: numerics}
	To numerically study the dynamical behavior of the traffic accident network model including the Hawkes process, we need to introduce a suitable implementation. For an overview we present Algorithm \ref{alg:Hawkes}.

 \begin{algorithm}[ht!]
    \caption{Traffic accident model algorithm using the Hawkes process}
    \begin{algorithmic}[1]
        \REQUIRE $T$, $\Delta t, \lambda_s,\lambda_d, \bar{d}$ and $ \rho_0^e, c_\text{road}^e, a_e,b_e$ for $e \in \mathcal{E}$
        \STATE $t=0$, $\rho = \rho_0$, $J=0$
        \WHILE{$t<T$}
        \STATE $u_1 \sim \mathcal{U}([0,1])$
            \IF{$u_1\leq \Delta t \lambda^*(t)$}   
                \STATE $J=J+1$
                \STATE $t_J= t$
                \STATE Decide on which road or junction the accident happens using $\mu_t^{\text{RI}}$ in \eqref{accmeasChoiceJunction}
                \IF{The accident happens on a road}
                    \STATE sample the accident position $p_J$ using $\mu_t^\text{pos,e}$ in \eqref{eq: mute}
                \ENDIF
                \STATE sample $c_J \sim \text{Beta}(2.62,3.53),s_J\sim Exp(\lambda_s)$ and $d_J\sim \bar{d}+ Exp(\lambda_d)$ 
            \ENDIF
            \STATE Adjust $c_a$ given $t_j,p_j,c_j,s_j,d_j,~j=1,\dots,J$
            \STATE Calculate one step of $\Delta t$ of the network traffic model using a Godunov scheme as presented in \eqref{eq:GodunovNetwork} and update $\rho^e, ~ e\in\mathcal{E}$
        \ENDWHILE
        \ENSURE The accident process $(t_j,p_j,c_j,s_j,d_j)_{j \in \lbrace 1,\dots,J\rbrace}$ and traffic densities $\rho^e, e\in \mathcal{E}$
    \end{algorithmic}
    \label{alg:Hawkes}
\end{algorithm}
	
	Starting with the Hawkes process, we state that the jump times of the Hawkes process can be obtained following the description in Section \ref{sec4}. For the measure $\mu_t^{RI}$ in \eqref{accmeasChoiceJunction} we apply a discrete inverse transformation method \cite{Asmussen2007} whereas for the position measure on a single road $\mu_t^\text{pos,e}$ in \eqref{eq: mute} a continuous inverse transformation method is reasonable. The same applies for the exponential distributions for the accident sizes and durations. For the Beta distribution, we exploit that a random variable $Z \sim Beta(\alpha',\beta')$ can be interpreted as $Z=\frac{X}{X+Y}$, where $X \sim \Gamma(\alpha',\theta), Y\sim \Gamma(\beta',\theta)$ for $\alpha',\beta',\theta >0$ which can be obtained out of an exponentially distributed random variable again (by an appropriate scaling). For further details of this approach see e.g. \cite{Asmussen2007,Georgii2009}. {In addition, \cite{Chen2016} further investigates a thinning algorithm and }the book by Laub et al. \cite{Laub2022} presents further common methods to sample the Hawkes process.\\
	The traffic dynamics is approximated on a time grid $(t^l)_{l \in \mathbb{N}},$ where $t^l = (l-1) \Delta t$ for some temporal stepsize $\Delta t>0$. 
	For the numerical discretization of the densities on the arcs we use a Godunov scheme combined with demand and supply functions for the fluxes at the vertices. The equidistant spatial grid on arc $e \in \mathcal{E}$ is given by $(x_{k,e})_{\lbrace k = 1, \dots, K_e \rbrace},$ where $x_{k,e} = (k + \frac{1}{2})\Delta x$ for some spatial stepsize $\Delta x>0$ and $x_{k,e}$ being the cell center. The left and right cell boundaries are given by $x_{k - \frac{1}{2},e} = k \Delta x$ and $x_{k + \frac{1}{2},e} = (k+1) \Delta x$, respectively. The stepsizes are chosen such that the CFL-condition
	\begin{align}
 \label{eq: CFLCondition}
		\Delta t \leq \frac{\Delta x}{\max_{e \in \mathcal{E}}|\sup_{x,\rho} \lbrace (F_e(x,\rho))_\rho\rbrace|}
	\end{align}
	holds true in every step for $(F_e)_\rho$ the partial derivative of $F_e$ with respect to the density variable.
	The Godunov scheme then reads for an arc $e \in \mathcal{E}$
	\begin{align}
 \label{eq:GodunovNetwork}
		\begin{split}
		    \rho_{1,e}^{l+1} &= \rho_{1,e}^l - \frac{\Delta t}{\Delta x} \left( G_e(\rho_{2,e}^l,\rho_{1,e}^l,x_{2,e},x_{1,e}) - F_e^\text{out}(a_e,\rho_{1,e}^l)  \right)\\
		\rho_{k,e}^{l+1} &= \rho_{k,e}^l - \frac{\Delta t}{\Delta x} \left( G_e(\rho_{k+1,e}^l,\rho_{k,e}^l,x_{k+1,e},x_{k,e}) - G_e(\rho_{k,e}^l,\rho_{k-1,e}^l,x_{k,e},x_{k-1,e}) \right)\\
		\rho_{K,e}^{l+1} &= \rho_{K,e}^l - \frac{\Delta t}{\Delta x} \left(F_e^\text{in}(b_e,\rho_{K,e}^l) - G_e(\rho_{K,e}^l,\rho_{K-1,e}^l,x_{K,e},x_{K-1,e}) \right),
		\end{split}
	\end{align}
	where the numerical flux is given by 
	\begin{align}
		\label{GGod}
		G_e(\rho_{k+1,e}^l,\rho_{k,e}^l,x_{k+1,e},x_{k,e}) = \min \left\lbrace F_e(x_{k+1,e}, \max\lbrace\rho^*,\rho_{k+1,e}^l\rbrace), F_e(x_{k_e}, \min\lbrace\rho^*,\rho_{k,e}^l\rbrace) \right \rbrace
	\end{align}
	with $F_e^\text{out}(\rho_{1,e}^l,a_e)$ being the outflow from the junction into arc $e$ and $F_e^\text{in}(\rho_{K,e}^l,b_e)$ the inflow into the junction from arc $e$ at time $t_l$. The initial data is determined by cell averages
	\begin{align*}
		\rho_{k,e}^0 = \frac{1}{\Delta x} \int_{x_{k-\frac{1}{2},e}}^{x_{k+\frac{1}{2},e}} \rho_{0,e}(x) dx.
	\end{align*}
	The function $G_e$ in (\ref{GGod}) has also an interpretation in terms of the demand and supply framework (see Section \ref{Sec: ExtNetwork}) for the cells on the road itself. 
	A single road can be considered as a sequence of cells having a certain capacity and transporting traffic flow. Thus, the demand for the available capacity in cell $k$ is given by $F_e(x_{k,e}, \min\lbrace\rho^*,\rho_{k,e}^l\rbrace),$ whereas the supply of the available capacity at the direct successor cell is given by $F_e(x_{k+1,e}, \max\lbrace\rho^*,\rho_{k+1,e}^l\rbrace)$ with $\rho^*$ denoting the flux maximizing density. The maximal flow in the demand and supply framework is then given by the minimum over both quantities. 
	In contrast to the theoretical investigations, we assume that the left- and rightmost roads in the network do not have infinite length. To buffer inflow into the traffic network, we install a queue whose evolution is given by the ordinary differential equation (ODE)
	\begin{align*}
		\frac{d}{dt} \hat{q}(t) = f^\text{in}(t) - F_{1}^\text{out}(a_1,\rho_1(a_1,t)),~~~ \hat{q}(0)=0,
	\end{align*}
	where $f^\text{in}(t)$ is the boundary inflow and $F_{1}^\text{out}(a_1,\rho_1(a_1,t))$ is the inflow into the traffic network at road 1. At time $t=0$ the queue is assumed to be empty. The ODE can be discretized using an explicit Euler scheme
	\begin{align*}
		\hat{q}^{l+1} =  \max \left\lbrace 0, \hat{q}^{l} + \Delta t f^\text{in}(t^l) - \Delta t F_1^\text{out}(a_1,\max \lbrace\rho^*, \rho_{1,1}^l \rbrace) \right\rbrace,
	\end{align*}
	choosing $\Delta t>0$ sufficiently small.
	The inflow into the queue $f^\text{in}$ is externally given, whereas the outflow again results in a demand and supply consideration restricted by 0 from below since the queue can never have a negative length.
	
	Next we present simulation results for the diamond network shown in Figure \ref{networkDiamond}. 
	
	\tikzset{
		pp/.style={
			rectangle,
			fill=orange!10,
			draw=orange
		}
	}
	\tikzset{
		ds/.style={
			ellipse,
			fill=blue!10,
			draw=blue
		}
	}
	\tikzset{
		cs/.style={
			rectangle,
			fill=green!10,
			draw=green
		}
	}
	\vspace{1cm}
	\begin{minipage}{0.58\textwidth}
		\begin{tikzpicture}[thick]
			
			\node[ds] (K1) {A};
			\node[pp,left = 0.1cm of K1] (K0) {queue};
			\node[ds,right = 0.8cm of K1] (K2) {B};
			\node[ds,above right = 1.2cm of K2] (K3) {C};
			\node[ds, below right = 1.2cm of K2] (K4) {D};
			\node[ds, below right = 1.2cm of K3] (K5) {E};
			\node[ds, right = 0.8cm of K5] (K6) {F};

			\path [->]
			(K1) edge[above] node{1} (K2)
			(K2) edge[above] node{2}(K3)
			(K3) edge[left] node{4}(K4)
			(K4) edge[above] node{6}(K5)
			(K3) edge[above] node{5}(K5)
			(K5) edge[above] node{7}(K6)
			(K2) edge[above] node{3}(K4);
		\end{tikzpicture}
		
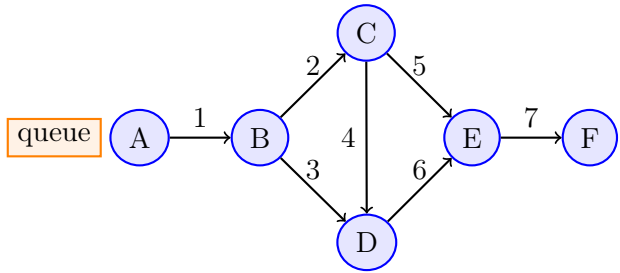
\captionof{figure}{The diamond traffic model.}
		
		\label{networkDiamond}
	\end{minipage}
	\begin{minipage}{0.38\textwidth}
		\centering
			\begin{tabular}{c|c|c}
				\hline
				arc & $c_\text{road}^e$ & $\rho_0^e$    \\
				\hline
				\hline
				1 & 0.7 & 0.4 \\
				\hline
				2 & 0.8 & 0.4\\
				\hline
				3 & 0.4 & 0.4\\
				\hline 
				4 & 0.5 & 0.8\\
				\hline 
				5 & 0.3 & 0.4\\
				\hline 
				6 & 0.8 & 0.8\\
				\hline 
				7 & 1 & 0.2\\
				\hline
			\end{tabular}
		\captionof{table}{Parameter choices for the diamond network.}
		\label{Parameters}
	\end{minipage}
	\vspace{0.2cm}\\
	
	It consists of six nodes (junctions) and seven directed arcs (roads) all of length 1 such that traffic flows from node $A$ to node $F$, where arc 4 only allows for flow from node C to node D. 
	A queue located upstream from node $A$ buffers potential inflow that cannot enter the system because it is already congested. Table \ref{Parameters} gives an overview about the parameters assumed for the network. The initial densities are constant. 
	The inflow is given by the function $f^\text{in}(t) = (0.13 + 0.052 \sin(t)) \mathbbm{1}_{[0,T-75]}(t)$ for $T=500$. 
	
	We denote by $\alpha_1=0.6,\alpha_2=0.5$ the distribution parameters at the junction B for the outflow to road 2 and at junction C to road 4, respectively. Directly, the shares floating out of junction B to road 3 is given by $1-\alpha_1$ and from C to road 5 by $1-\alpha_2$. The priority rightway parameters at the junctions D and E are denoted by $q_1=0.5$ from road 3 ($1-q_1$ from road 4) and $q_2=0.4$ from road 5 ($1-q_2$ from road 6). We choose $\Delta t = \frac{1}{100} = \Delta x$ and the velocity function to the common choice of $v(\rho) = 1- \rho$. The CFL-condition in \eqref{eq: CFLCondition} is satisfied since the road capacities in Table \ref{Parameters} do not exceed the value of 1 and therefore $|(F_e(x,\rho))_\rho|\leq 1$.
	For the accident generation we choose $\gamma = 0.5, \gamma_v = 0.2$ scaling the accident risk on the roads and at the junctions and $\alpha = 0.1 , \beta = 2, \tilde{\beta} = 24, \nu =0$ for the parameters of the self-excitation property. An exponential distribution with parameter {$\lambda_s=20$} is chosen for the accident size. For the duration of an accident we set $d_{\text{add}} \sim Exp(\frac{1}{2})$ and $\bar{d}=1$. These two choices deviate at a first glance form the choices in Section \ref{sec4}, but can be obtained when scaling the road properly. For the capacity reduction we stick to the suggested choice for the beta distribution from the data analysis in Section \ref{sec4} and set $c \sim Beta(2.66,3.53)$. 
	
	
	\subsection{Accident analysis}
\begin{figure}[ht!]
		\centering
		\includegraphics[scale = 0.74]{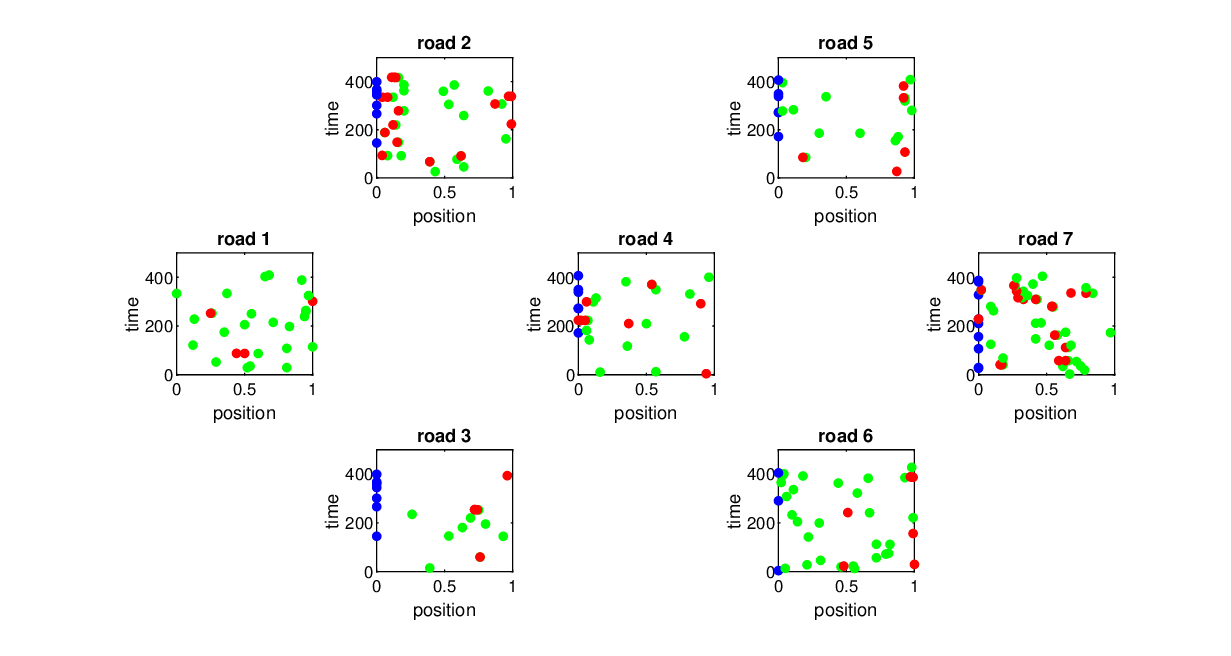}
		\caption{The positions of accidents in time and space for the diamond network visualized by green (background accidents), red (self-excitation accidents) and blue (accidents at junctions) circles.}
		\label{fig:accanalysis1}
	\end{figure}

	A visualization of the accident behaviour in the diamond network is presented in Figure \ref{fig:accanalysis1}, {where the roads are numbered as in Figure \ref{networkDiamond}}. For each road of the diamond network the circles show accident positions in {space (horizontal axis) and time (vertical axis)}. The green circles reflect background accidents due to the flux-share of the Hawkes process whereas the red circles show self-excitation accidents. The blue circles describe an accident at the junction and are always depicted at the beginning of the adjacent outgoing roads. For junctions with more than one outgoing road (e.g. at junction B with road 2 and road 3) the accidents are displayed at both roads at the same time points.
	
	We observe that while the green circles do not follow a particular pattern, the red circles often accumulate close to other green or red accidents. They are always positioned to the top left of a background (primary) accident since a secondary self-excitation accident can only happen behind a primary accident in space and only after a primary accident in time. Furthermore, we identify that overall on road 1 and 7 there are more accidents than e.g. on road 4 and 5, which is mainly due to the larger traffic flow since all of the vehicles have to pass road 1 and 7, whereas road 5 can be bypassed using road 6. On road 6 close to position $x=1$ and very small times we see that there are self-excitation accidents which are now caused by a primary accident on road 7, which shows that the accident risk laps over the different roads. 
	
\begin{figure}[ht!]
		\begin{minipage}{0.49\textwidth}
			\centering
			\includegraphics[scale=0.54]{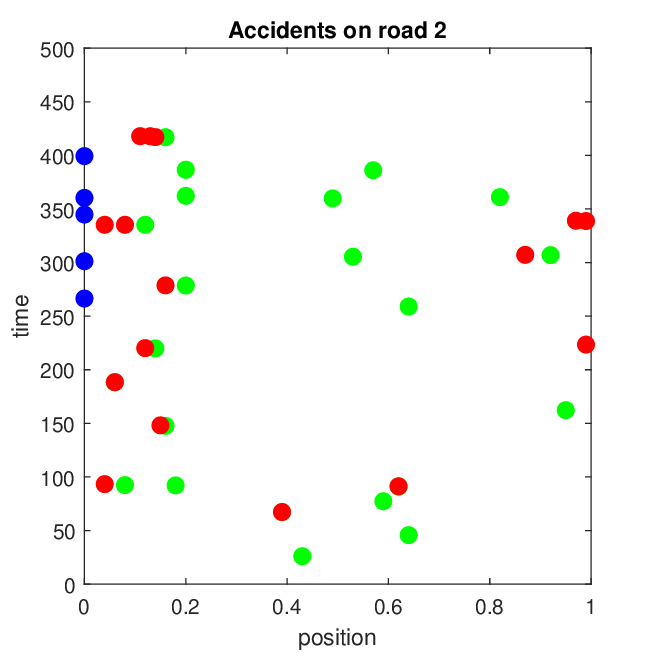}
		\end{minipage}
		\begin{minipage}{0.49\textwidth}
			\includegraphics[scale = 0.62]{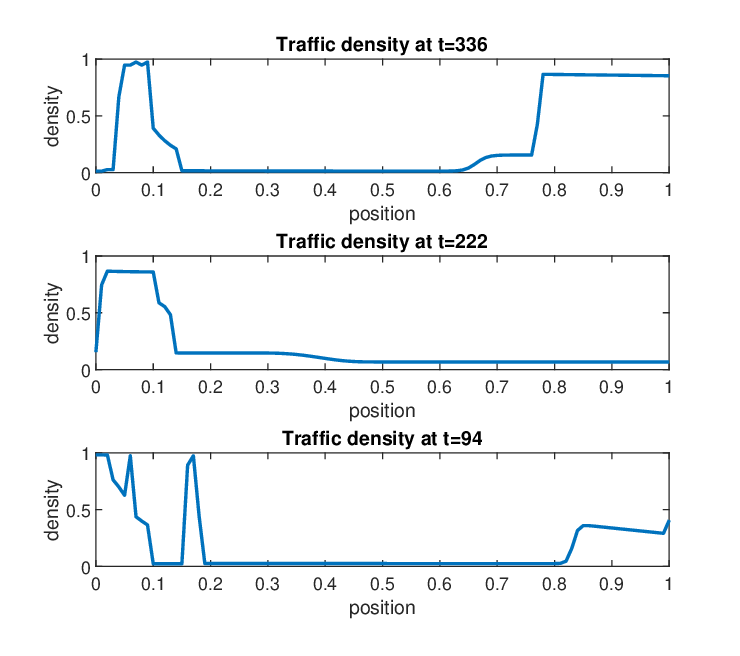}
		\end{minipage}
		\caption{Accidents on road 2 on the left, cf. Figure~\ref{fig:accanalysis1},
			and corresponding traffic densities on road 2 for $t=94,222,336$ on the right.}
		\label{fig:accanalysis2}
	\end{figure}

	
	To relate accidents and densities we focus on road 2 and compare the accident occurrence with the traffic densities in Figure \ref{fig:accanalysis2}. Around $t=94$ and $t=336$, we observe three accidents in the road section $[0,0.2]$ which lead to three increases in the density there. Since only very few vehicles are able to pass this area, the density is close to zero on $[0.2,0.8]$ and $[0.2,0.6]$. For $t=222$ around $x=0.15$ a pair of a background and self-excitation accident decreases the road capacity which again leads to congestion and a high density regime at the beginning of the road.
	
 {
 This study underlines the effect of self-excitation in the simulation that originally was indicated in the histogram of Figure \ref{intTimes1}. 
 
 In a next step, we validate the construction of the accident times using a comparison to real data. The idea is to compare the relative number of accidents that fall into hourly intervals within one day. Since accidents happen rarely, we have to choose a sufficiently large region for our numerical analysis. Therefore, we stick to the one-year \textit{accident data} set of Great Britain\footnotemark[1]\footnotetext[1]{\url{https://data.gov.uk/dataset/cb7ae6f0-
			4be6-4935-9277-47e5ce24a11f/road-safety-data},   \\ Accessed: 2024-11-05} already used in Section \ref{sec3}.
 Since a simulation of the entire traffic evolution of one year in Great Britain is computationally not feasible, we use the diamond network as an exemplary traffic system and feed it with hourly changing inflow. The inflow is also based on real data information. Unfortunately, the data quality of corresponding British vehicle counting data was not good enough, thus the inflow data is derived from \textit{vehicle count data} in Germany\footnotemark[2]\footnotetext[2]{\url{https://www.bast.de/DE/Verkehrstechnik/Fachthemen/v2-verkehrszaehlung/Verkehrszaehlung.html?nn=1817946}, Accessed: 2024-11-05}. We averaged data from two highway vehicle counting stations in Appenweier and Hockenheim as well as traffic data from two federal roads in Rastatt and Breisach to take different types of roads into consideration. 
 Generally, we stick to the parameters of the diamond network introduced before, but interpret one time unit as one hour and therefore increased the road capacities from Table \ref{Parameters} by factor 10, set $T=24$ and reduced the accident duration parameter to $\bar{d}=0.5, \lambda_d = 1$ accordingly. For the simulation of the traffic accident model we simulate 500 days for weekday traffic inflow data and also 500 days for Sunday inflow data. Figure \ref{fig: weekdayDataTime} presents the averaged relative number of accidents falling into the 24 hourly intervals for the traffic accident model using the Hawkes process to determine the accident times in blue and the real data in the red dots. The same comparison for Sundays can be found in Figure \ref{fig: sundayDataTime}. We may observe the different properties of traffic accidents on Monday to Friday and Sunday, where on a weekday there is a two-peak structure, one for the rush hour in the morning around 8am and one in the afternoon around 5pm, whereas on Sunday the highest accident risk is reached in the afternoon. The model data captures well the main characteristics of the real data, decreases to very low numbers at night and shows the highest shares around the rush hour times. Only the most extreme outliers are not captured perfectly. On Sundays, there is fewer traffic in the morning and the model still predicts correctly the most accidents in the afternoon hours. It merely underestimates the number of accidents in the night hours. \\ 
 \begin{minipage}[ht!]{0.48\textwidth}
 \centering
 \includegraphics[scale=0.52]{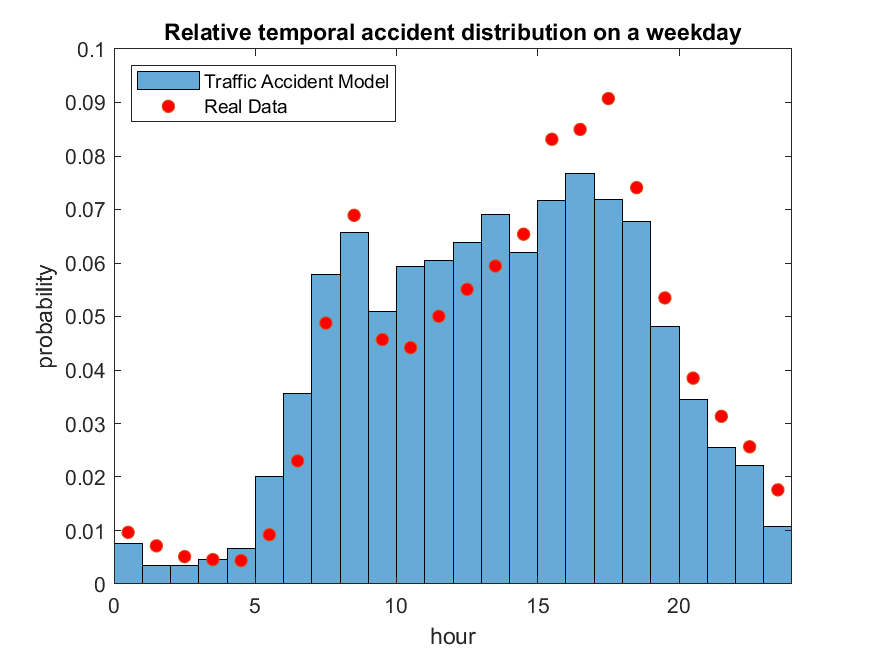}
 \captionof{figure}{Weekday accident time comparison of the accident traffic model and real data.}
     \label{fig: weekdayDataTime}
 \end{minipage} ~~\begin{minipage}{0.48\textwidth}
 \centering
 \includegraphics[scale=0.52]{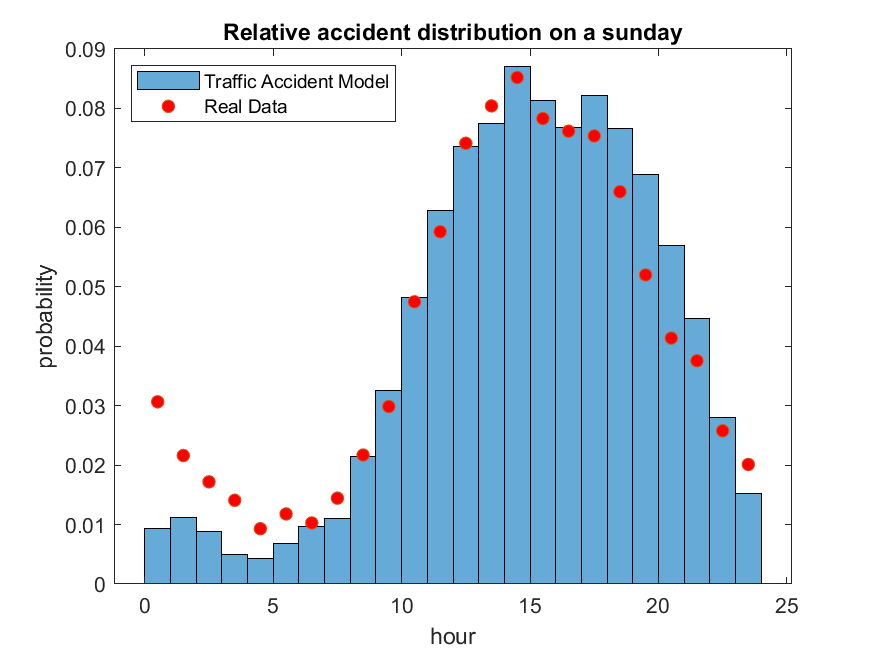}
 \captionof{figure}{Sunday accident time comparison of the accident traffic model and real data.}
     \label{fig: sundayDataTime}
 \end{minipage}
 
 }

	\subsection{Risk Measures}
	\label{Sec: riskMeasures}
	
	In the following we analyze the network performance dependent on the choices of the distribution parameters $\alpha_1$ and $\alpha_2$ by introducing risk measures. These are evaluated using the parameters from before but considering $\gamma=0.1, \gamma_v = 0.04$ to reduce the computational effort and a Monte Carlo simulation of 1000 runs to capture different accident scenarios with a time horizon of $T=150$. The use of risk measures has been previously presented e.g. in \cite{Friedrich2022,Reilly2015,Treiber2013}.
	
	\subsubsection{Total travel time}
	First, we discuss the \textbf{total travel time (TTT)} which measures how much time the vehicles spend in the system. The measure is composed of the time vehicles drive in the network and wait until they can enter the network in the queue. 
	\begin{align*}
		TTT = \sum_{e = 1}^7 \int_0^1 \int_0^T \rho_e(x,t)dt dx + \int_0^T \int_0^T q(r)dr ds
	\end{align*}
	As before $q$ denotes the length of the queue prior to road 1 and $\rho_e(x,t)$ is the density of road $e$ at $x$ and time $t$. 
	Generally, lower total travel times are desirable because less time is required to get through the traffic network which additionally reduces for example traffic noise and emissions. Note that the total travel time can only be compared if we ensure that the system is empty at all terminal times. 
	
	We fix the rightway parameters to $q_1 = 0.5, q_2 = 0.4$ and vary the distribution parameters in steps of $0.1$ between 0 and 1 to compare the averaged total travel times of the Monte Carlo simulation in Figure \ref{fig11}. We observe the lowest values of the total travel time for choices of $\alpha_1 \in [0.5,0.6]$ and $\alpha_2\in [0.3,0.4]$. This means that at junction B a slightly larger share of cars is governed to road 2 and at junction C more vehicles move to road 5. Comparing the averaged accident total travel times to a benchmark setting with the same parameters but no accidents, Figure \ref{fig11bm} shows that on the one hand for all parameter choices travel times increase with the introduction of accidents, where e.g. for the minimal values we observe an increase of 50 percent in the total travel times. The largest increases are observed for configurations with high values in $\alpha_1$ and $\alpha_2$ which is consistent with the results in the following Section \ref{sec: Noa}. On the other hand, the parameter configuration leading to the minimum is slightly shifted and attained now for smaller $\alpha_1$ little and significantly larger $\alpha_2$. Generally, there is a wider range of parameter configurations that lead to comparably low travel times in the setting without accidents.\\
	\begin{minipage}[ht!]{0.48\textwidth}
		\centering
		\includegraphics[scale=0.52]{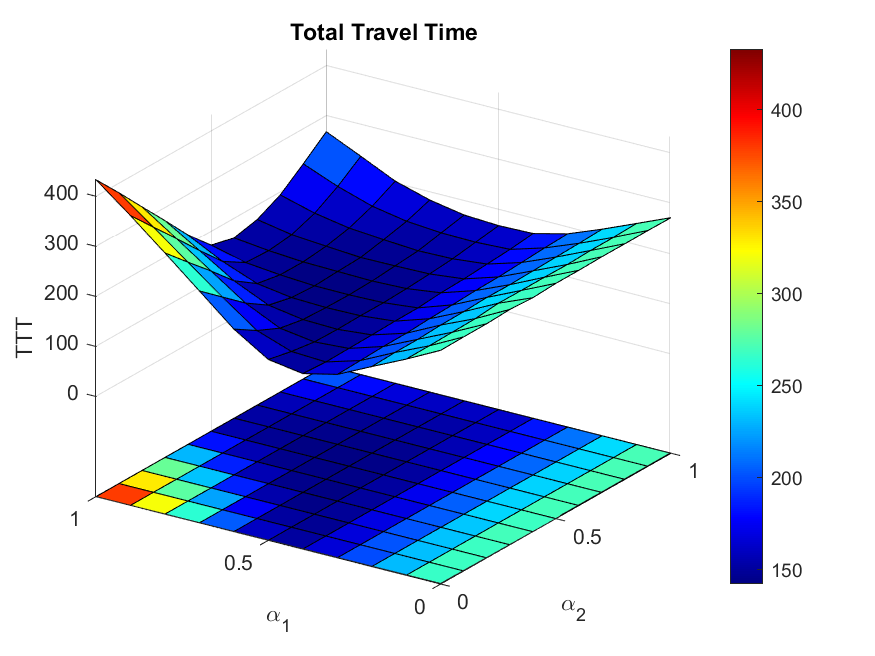}
		\captionof{figure}{Total travel time dependent on the parameter choices of $\alpha_1$ and $\alpha_2$ with accidents.}
		\label{fig11}
	\end{minipage}~
	\begin{minipage}[ht!]{0.48\textwidth}
		\centering
		\includegraphics[scale=0.52]{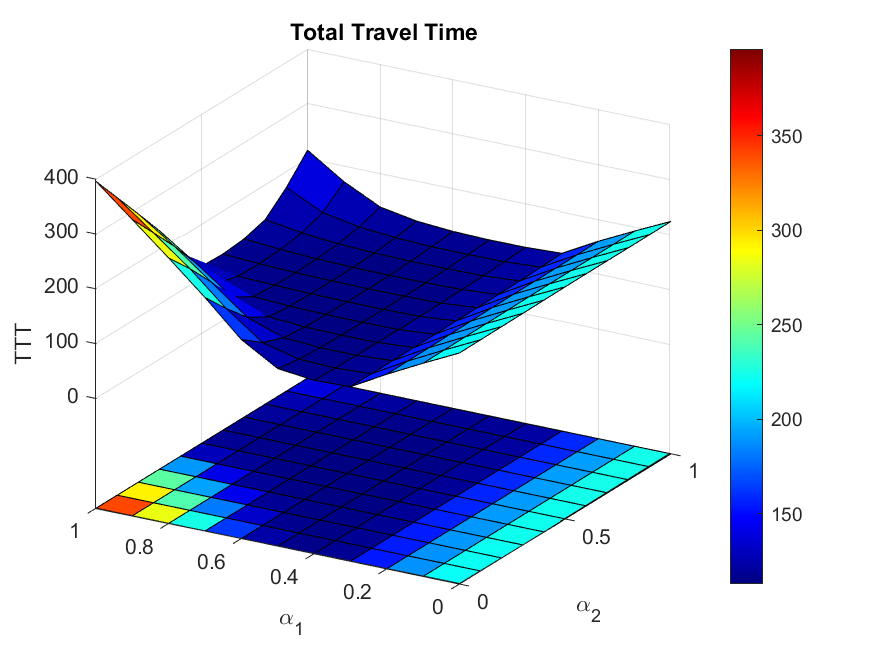}
		\captionof{figure}{Total travel time dependent on the parameter choices of $\alpha_1$ and $\alpha_2$ in a benchmark setting without accidents.}
		\label{fig11bm}
	\end{minipage}

	\subsubsection{Number of accidents}
	\label{sec: Noa}
	In this subsection, we take a look on the road safety and compare the total number of accidents. As before, we leave the rightway parameters fixed to $q_1 = 0.5, q_2 = 0.4$ and consider the number of accidents dependent on choices of $\alpha_1$ and $\alpha_2$. In contrast to the total travel time results, Figure \ref{fig51} shows a monotonicity in both arguments. The largest number of accidents is obtained for $\alpha_1=\alpha_2 = 1$ which is the case in which all vehicles choose the route via the roads 2,4 and 6 to get from junction B to junction E. The number of accidents decreases for lowering both parameters and reaches the lowest values for $\alpha_1, \alpha_2$ close to 0. Even though, this leads to low numbers of accidents, we remark that the total travel times are high in these scenarios indicating that road safety and network efficiency are not necessarily positively correlated. \\
	\begin{center}
		\begin{minipage}{0.48\textwidth}
			\centering
			\includegraphics[scale=0.5]{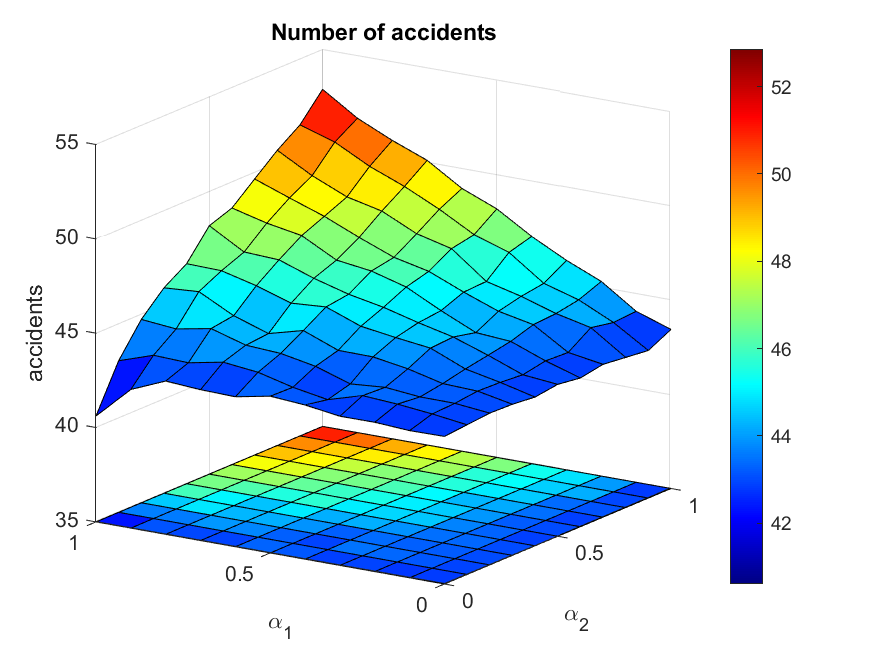}
			\captionof{figure}{The number of accidents for different choices of $\alpha_1$ and $\alpha_2$.}
			\label{fig51}
		\end{minipage}
		%
		\hspace{0.1cm}
		
	\end{center}

	\subsubsection{Probability of an empty system at some time $t$}
\begin{figure}[ht!]
		\centering
		\includegraphics[scale=0.57]{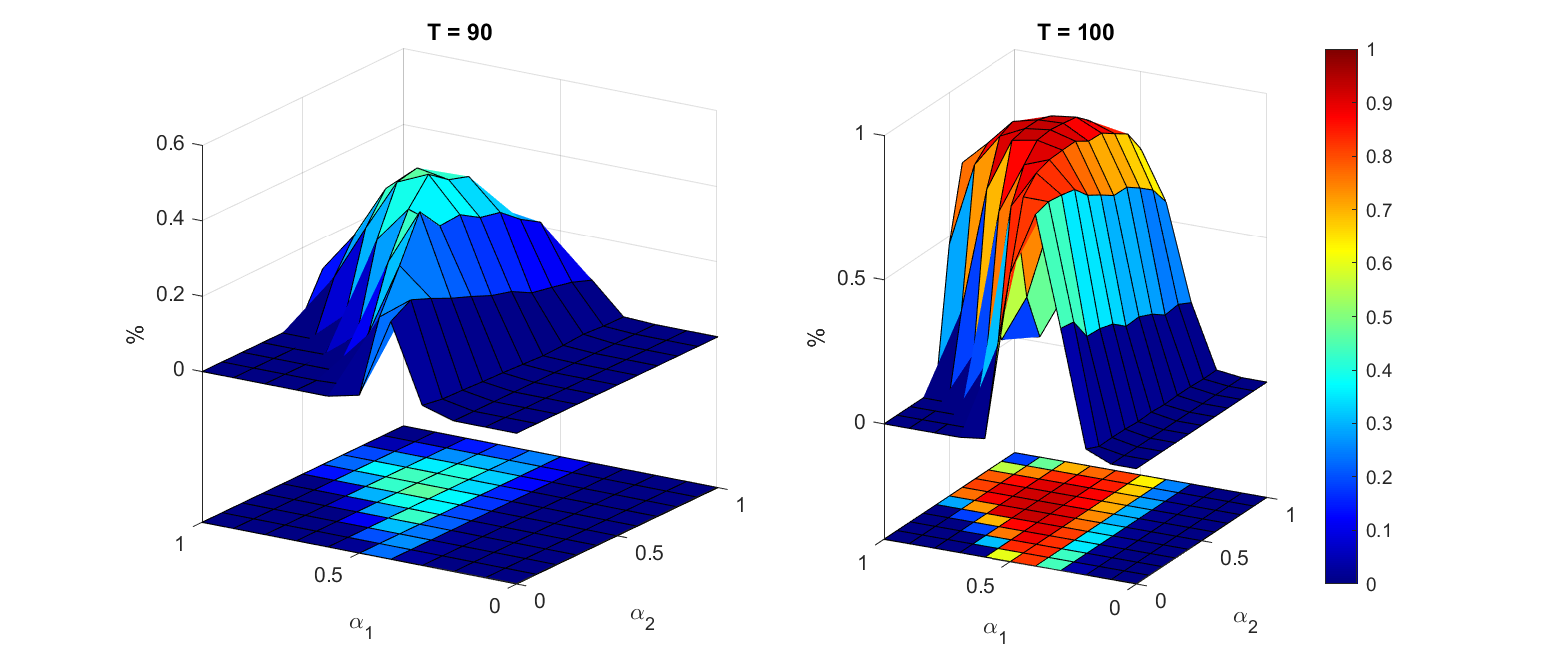}
		\caption{Probabilities of an empty system dependent on the parameter choices of $\alpha_1$ and $\alpha_2$ for $T\in \lbrace 90,100 \rbrace$.}
		\label{fig41}
	\end{figure}

\begin{figure}[ht!]
		\centering
		\includegraphics[scale=0.57]{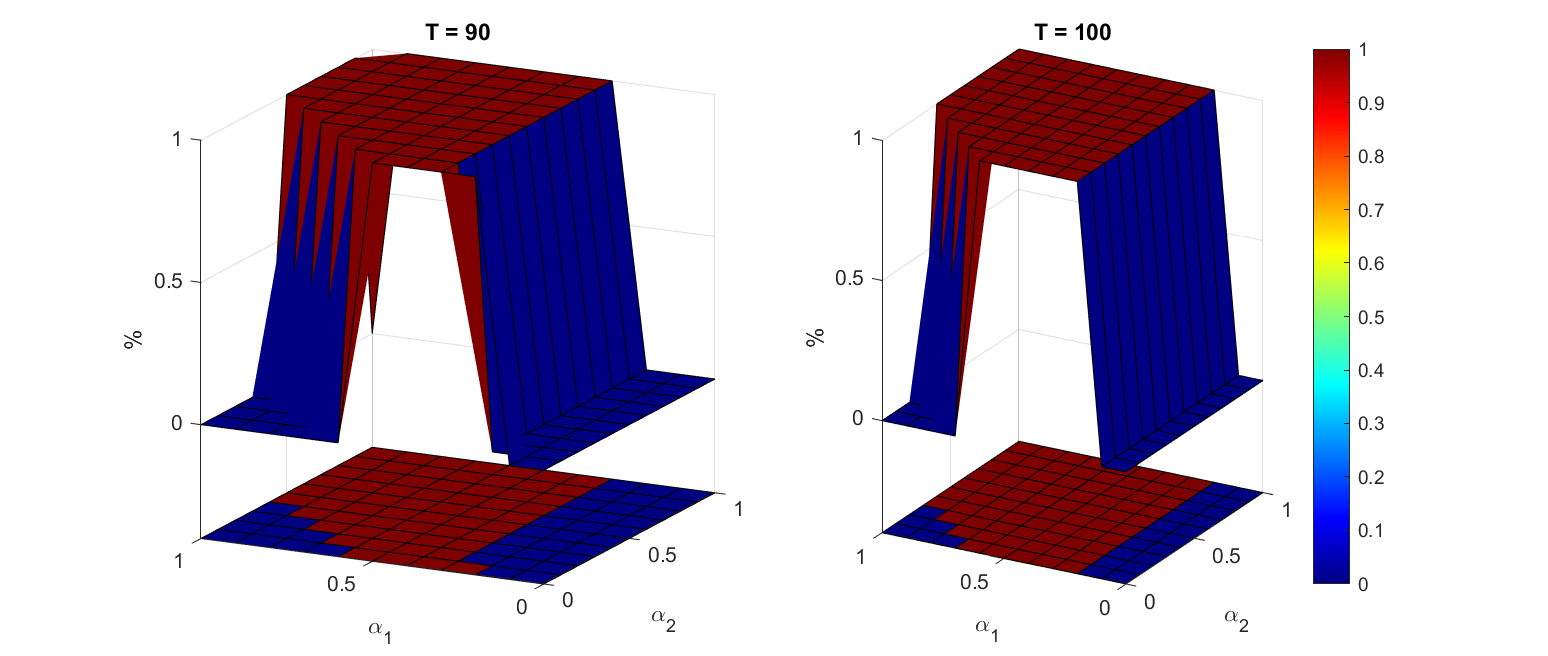}
		\caption{Probabilities of an empty system dependent on the parameter choices of $\alpha_1$ and $\alpha_2$ for $T\in \lbrace 90,100 \rbrace$ in a (deterministic) benchmark framework without accidents.}
		\label{fig411}
	\end{figure}

	Since there is no more inflow into the system after $t=75$, it runs out of vehicles after some time. This particular time is called time of an empty system (ToES).
	Here, we focus not only on averaged values but on the particular distributions of the times of an empty system and determine the probability that for a particular time $t>75$ the system is empty, i.e. $P(ToES\leq t)$. Note that the delays until the system is empty are mainly caused by arising accidents. 
	Keeping the rightway parameters fixed to $q_1 = 0.5, q_2 = 0.4$ we consider the probabilities for an empty system at $t \in \lbrace 90,100 \rbrace$ depending on the choices of $\alpha_1, \alpha_2$ in Figure \ref{fig41}. The probability of an empty system is approximated by a Monte Carlo simulation again. With increased times also the probabilities of an empty system naturally increase. For $T=90$ the largest probability for an empty system are attained at $\alpha_1=0.5, \alpha_2 = 0.6$ whereas for $T=100$ we observe a plateau with a probability of more than 90 percent of an empty system for $\alpha_1 \in [0.5,0.7]$ and $\alpha_2 \in [0.4,0.8]$. Low values of $\alpha_1$ seem to lead to a very congested setting in which the system is almost never empty even after $t=100$. The same holds true for combinations of $\alpha_1$ close to 1 and $\alpha_2$ close to 0. 
	We also present a reference framework without accidents in Figure \ref{fig411} which due to its deterministic structure only allows for probabilities of 0 or 1. At least for $T=100$ the most configurations lead to an empty system. Comparing the results to the accident framework, the impact of the accidents and the corresponding congestion get visible by the significantly lower shares of an empty system. 

	\subsection{Rerouting strategies}
	We consider again the diamond network from Figure \ref{networkDiamond}. There are basically three ways to get through the network, i.e., route 1 via the roads 1,2,5 and 7, route 2 via the roads 1,3,6,7 and route 3 using the inner connection over the roads 1,2,4,6,7. The routing preferences of the individual cars are represented by the distribution parameters $\alpha_1$ and $\alpha_2$. Let us assume $\alpha_1 = 0.65$ and $\alpha_2 = 0.3$. This corresponds to $45.5$ percent of the driver taking route 1, $35$ taking route 2 and $19.5$ percent choosing route 3. Additionally, we assume that $\frac{4}{7}$ of the drivers choosing route 1 are considered to be flexible drivers, meaning that in case of large congestion on road 5 and a free traffic situation on the roads 4 and 6, they change their plan and choose route 3 instead of route 1. This can be considered as the second distribution parameter at junction C changing to $\alpha_2' = 0.7$. After the congestion or the serious accident on road 5 is resolved, the distribution parameter at junction C changes back to $\alpha_2 = 0.3$. In practice, digital road signs at junctions can announce congestion on road 5 and recommend a detour via the roads 4 and 6. 
	
	In our particular case a detour is recommended in two different cases: large congestion on road 5 or a serious accident on road 5. We measure congestion using the following measure $CM_e(t)$ for road $e$
	\begin{align*}
 \vspace{-0.5cm}
		CM_e(t) =\max \left\lbrace  \int_0^1 \rho^e(x,t)-  \frac{f^e(\rho^e(x,t),x,t)}{v_\text{ref}}dx,0\right\rbrace,
	\end{align*}
	where we choose $v_{\text{ref}} = \frac{v_\text{max}^e}{2} = \frac{1}{2}$. The road is considered to be congested if $CM_e(t)$ exceeds the value 0.25.
	A serious accident is defined to be an accident for which the capacity reduction $c$ is larger than 0.8. But since it is not useful to recommend a detour in one of these cases if road 4 or 6 is also congested or has a serious accident, the detour can only be recommended if $\max \lbrace CM_4(t), CM_6(t) \rbrace \leq 0.25$ and there is no accident on road 4 and 6 having a capacity reduction of larger than 0.8. Apart from the scenario with some flexible drivers (scenario II, $\alpha_2' = 0.7$), we also consider a situation (III) with all drivers reaching junction C are flexible drivers ($\alpha_2' = 1$). 
 Since in this framework the distribution parameters additionally depend on the traffic density, the existence of solutions cannot be ensured by the theoretical results used before, cf. also \cite{Garavello2016, Garavello2006}.
 
 Comparing the introduced scenarios with the standard benchmark setting without any detour recommendation in a Monte Carlo simulation with 2000 runs we use the risk measures introduced before and show their values in Table \ref{detour1} and \ref{detour2}. Scenario I describes the standard case where we do not consider any flexible drivers. If no drivers adapt their behaviour the total travel time amounts to 256.45. We observe that we can reduce this value having some drivers willing to take the detour in case of a congestion on road 5 of almost 5 percent $(\alpha_2'=0.7)$ or even more in case of $\alpha_2'=1$. To show that the choice of $\alpha_2=0.7$ is not superior in general, we simulate scenario IV and obtain an even longer total travel time. But not only the total travel time benefits from intelligent routing, also the probabilities for an empty system at $T \in \lbrace 90,100,110 \rbrace$. The values for the detour frameworks are always superior compared to the one without flexible drivers. 
 Hence considering the numbers of accidents in Table \ref{detour2} there is no clear conclusion. As expected, we observe a decrease in accidents on road 5 if more vehicles take the detour via road 4 and 6, which then leads to more accidents on these two roads. Overall, detours might slightly increase the number of accidents but also reduce the total travel time. 
	\begin{table}[ht!]
		\centering
			\begin{tabular}{r|c||c||c|c|c}
				\hline
				&& TTT & $P(ToES<90)$ & $P(ToES<100)$ & $P(ToES<110)$ \\
				\hline
				\hline
				I&$\alpha_2 = \alpha_2' = 0.3$ & 256.45 & 0.004 & 0.079 & 0.326\\
				\hline
				\hline
				II&$\alpha_2 = 0.3, \alpha_2' = 0.7$  & 246.05 & 0.009 & 0.152 & 0.444 \\
				\hline
				III&$\alpha_2 = 0.3, \alpha_2' = 1$ & 243.14 & 0.011 & 0.189 & 0.469\\
				\hline
				\hline
				IV&$\alpha_2 = \alpha_2'=0.7$ & 263.04 & 0.012 & 0.143 & 0.396 
			\end{tabular}
		\captionof{table}{Total travel time and the probabilities of an empty system in the detour scenarios for junction C.}
		\label{detour1}
	\end{table}
	
	\begin{table}[ht!]
		\centering
			\begin{tabular}{r|c||c|c|c|c|c|c|c}
				\hline
				& & road 1 & road 2 & road 3 & road 4 & road 5 & road 6 & road 7   \\
				\hline
				\hline
				I&$\alpha_2 = \alpha_2' = 0.3$ & 8.09 & 5.62 & 3.09 & 2.14 & 4.33 & 5.79 & 10.39\\
				\hline
				\hline
				II&$\alpha_2 = 0.3, \alpha_2' = 0.7$  & 8.02 & 5.64 & 3.07 & 2.50 & 4.01 & 6.21 & 10.47\\
				\hline
				III&$\alpha_2 = 0.3, \alpha_2' = 1$ & 8.08 & 5.55 & 3.07 & 2.62 & 3.85 & 6.28 & 10.48 \\
				\hline
				\hline
				IV&$\alpha_2 = \alpha_2'=0.7$ & 8.23 & 5.69 & 3.18 & 4.36 & 2.18 & 8.09 & 10.57
			\end{tabular}
		\captionof{table}{Averaged numbers of accidents on the roads in the detour scenarios at junction C.}
		\label{detour2}
	\end{table}

	\section{Conclusion}
	This work introduces a traffic network model that captures traffic accidents and their self-excitation property. Accidents are reflected by capacity reductions in the model, whereas the traffic flow influences the background accident risk such that both ingredients are related to each other. A large emphasis is put on the data validation of the self-excitation property and all remaining accident parameters. The numerical framework presents the influences of traffic accidents on risk measures such as the total travel time or the number of accidents varying the distribution parameters of the network. A first insight in rerouting strategies is given which allows for further discussions on the \textit{best} choices of these parameters to reduce social costs. Therefore, future research might include the consideration of multi-objective optimization problems.
	
	\section{Acknowledgment}
	Simone G\"ottlich was supported by the German Research Foundation (DFG) under grant GO 1920/11-1 and 12-1.
	
	\section{Conflict of interest}
	The authors declare no competing interests.
	
	
	
	\printbibliography
\end{document}